\documentclass[12pt]{article}
\usepackage{amsmath}
\usepackage{graphicx,psfrag,epsf}
\usepackage{enumerate}
\usepackage{natbib}
\usepackage{url} 
\usepackage{amsfonts}
\usepackage{amsthm,amssymb}
\usepackage{verbatim}
\usepackage{bbm}

\usepackage{color}

\newcommand{\blind}{0}
\newcommand{\bo}[1]{\boldsymbol{#1}}

\newcommand{\E}{\mbox{E}}

\newcommand{\var}{\mbox{var}}
\newcommand{\cov}{\mbox{cov}}
\newtheorem{theorem}{Theorem}
\newtheorem{definition}{Definition}
\newtheorem{corollary}{Corollary}

\newcommand{\bmat}{\begin{pmatrix}}
\newcommand{\emat}{\end{pmatrix}}

\begin{document}

\def\spacingset#1{\renewcommand{\baselinestretch}%
{#1}\small\normalsize} \spacingset{1}


\if0\blind
{
  \title{\bf Extracting conditionally heteroscedastic components using ICA}
  \author{Jari Miettinen \\
		Department of Signal Processing and Acoustics, Aalto University, Espoo, Finland\\
    and \\
    Markus Matilainen \\
    Department of Mathematics and Statistics, University of Turku, Finland\\
    and \\
    Klaus Nordhausen \\
    Institute of Statistics \& Mathematical Methods in Economics, \\
Vienna University of Technology, Austria\\
    and \\
    Sara Taskinen \\
    Department of Mathematics and Statistics, University of Jyvaskyla, Finland}
  \maketitle
} \fi

\if1\blind
{
  \bigskip
  \bigskip
  \bigskip
  \begin{center}
    {\LARGE\bf Extracting conditionally heteroscedastic components using ICA}
\end{center}
  \medskip
} \fi

\bigskip

\begin{abstract}
	In the independent component model, the multivariate data is assumed to be a mixture of mutually independent latent components, and in independent component analysis (ICA) the aim is to estimate these latent components. In this paper we study an ICA method which combines the use of linear and quadratic autocorrelations in order to enable efficient estimation of various kinds of stationary time series. Statistical properties of the estimator are studied by finding its limiting distribution under general conditions, and the asymptotic variances are derived in the case of ARMA-GARCH model. We use the asymptotic results and a finite sample simulation study to compare different choices of a weight coefficient.  As it is often of interest to identify all those components which exhibit stochastic volatility features we also suggest a test statistic for this problem. We also show that a slightly modified version of principal volatility components (PVC) can be seen as an ICA method. Finally, we apply the estimators in analyzing a data set which consists of time series of exchange rates of seven currencies to US dollar. Supplementary material including proofs of the theorems is available online.
\end{abstract}
\noindent%
{\it Keywords:}  ARMA-GARCH Processes, Asymptotic Normality, Autocorrelation, Blind Source Separation, Minimum Distance Index, Principal Volatility Component
\vfill

\newpage
\spacingset{1.45} 
\section{Introduction}
\label{sec:intro}

In independent component analysis (ICA), the goal is to find a linear transformation of the multivariate data set which has mutually independent components. The purpose of carrying out ICA can be to separate interesting components from the noise, or to shift from a multivariate analysis to multiple univariate analyses. The latter is of particular interest when analyzing multivariate time series. As an example consider multivariate GARCH models, which according to \cite{Bollerslevetal1994} and \cite{Chibetal2006}, tend to have an overwhelming number of parameters, or be unreasonably simplified. In such case, it may be beneficial to reduce the dimension of data before carrying out further analyses. This was also suggested in~\cite{HuTsay2014}, where a method called {\it principal volatility component} (PVC) analysis was introduced.  For approaches applying independent component analysis in the context of multivariate GARCH modelling, see~\cite{Wu2005},~\cite{MattesonTsay2011} and~\cite{Hai2017}, for example.

Most of the ICA methods rely on making the marginal densities of the components maximally non-Gaussian. In the case of times series data, it is natural to make use of the temporal dependence. Arguably, the most famous one of such methods is SOBI (second-order blind identification) \citep{Belouchranietal1997}, which uses approximate joint diagonalization of a set of autocovariance matrices to find the transformation into independent components. SOBI performs well when the independent components have nonzero linear autocorrelations, but it fails to utilize volatility clustering. On the other hand, ICA estimators which are tailored for time series with volatility clustering, e.g. \citet{Hyvarinen2001,Shietal2009,Matilainenetal2015,Matilainenetal2017}, do not utilize information coming from linear autocorrelations to their full extent. In this paper we introduce an ICA method which we denote generalized SOBI (gSOBI). The method uses both linear and quadratic autocorrelations in order to find various kinds of time series. Furthermore, following \citet{HuTsay2014} who argued that especially in multivariate econometric times series the main interest is to identify those latent components which do not contain volatility clustering, we develop tests to identify such corresponding independent components. In this context we also show how PVC can be used to solve the ICA problem and suggest a small modification to make the method affine equivariant.

The paper is organized as follows. In Section~\ref{sec:ICA} we define the independent component model and discuss some ICA methods. In Section~\ref{sec:ARMAGARCH} we recall the ARMA-GARCH model which covers a wide collection of stationary time series. Section~\ref{sec:GSOBI} includes the formal definition of the generalized SOBI estimator, asymptotic results of the estimator, and tests for linear autocorrelation and volatility clustering. In Section~\ref{sec:pvca} we show that PVC can also be categorized as an ICA method utilizing volatility clustering. Finite sample properties of the gSOBI estimator, the PVC estimator, and the tests are studied in Section~\ref{sec:simulations}. In Section~\ref{sec:data}, the methods are applied to a dataset of exchange rates of seven currencies to US dollar, and the article is concluded in Section~\ref{sec:conc}.

\section{Independent component analysis}
\label{sec:ICA}

In the independent component model
\[
\bo x_t=\bo\Omega\bo s_t, \ \ \ t=0,\pm 1,\pm 2,\dots
\]
we denote the observable $p$-variate time series by $\bo x:=(\bo x_t)_{t=0,\pm 1,\pm 2,\dots}$, the full-rank  $p\times p$ mixing matrix by $\bo\Omega$, and the latent $p$-variate time series by $\bo s:=(\bo s_t)_{t=0,\pm 1,\pm 2,\dots}$. Given a realization $\bo x_1,\dots,\bo x_n$ of the process $(\bo x_t)_{t=0,\pm 1,\pm 2,\dots}$, the target in independent component analysis is to estimate the mixing matrix, or the unmixing matrix $\bo\Gamma=\bo\Omega^{-1}$ which transforms $\bo x$ back to $\bo s$ with $\bo s_t=\bo\Gamma\bo x_t$. Under certain conditions on $\bo s$, the unmixing matrix is indeed identifiable up to scales, signs and order of its rows. The key assumption is that the components of $\bo s$ are mutually independent. In addition, at most one of the components can be iid Gaussian, and due to the scale ambiguity, we assume for simplicity that the components of $\bo s$ have unit variances.

The usual strategy in ICA is first to standardize the data,
\[
\bo x_t^{st} = \bo S^{-1/2}(\bo x)(\bo x_t-\E(\bo x_t)),
\]
where $\bo S(\bo x)$ denotes the covariance matrix functional of $\bo x$. For the standardized data we then have that
\[
\bo s_t=\bo U\bo x_t^{st}
\]
for some orthogonal $p\times p$ matrix $\bo U$, see \cite{Miettinenetal2015}.

The ICA literature comprises a vast amount of estimators with different strategies to find the independent components. For the present, methods which only use marginal densities of the components and ignore the temporal or spatial order of the components, have been the most popular choices in ICA applications. However, most often the data in the ICA applications are time series or otherwise structured, and it is well-known that in these cases one can achieve better performance with methods which utilize the important data properties.

The standard way to inspect temporal dependence in a univariate time series $(x_t)_{t=0,\pm 1,\pm 2,\dots}$ is to compute autocovariances
\[
\E\left((x_t-\E(x_t))(x_{t+\tau}-\E(x_{t+\tau})\right),\ \tau>0.
\]
If the autocovariances of the independent components are nonzero, we can use joint diagonalization of autocovariance matrices in ICA \citep{Tongetal1990,Belouchranietal1997}. Notice that zero autocovariances do not imply the absence of temporal dependence, as, for example, economic and financial time series often exhibit volatility clustering, i.e., periods of lower and higher volatility. While the linear autocovariances might be all zero, the quadratic autocovariances
\[
\E\left((x_t-\E(x_t))^2(x_{t+\tau}-\E(x_{t+\tau}))^2\right),\ \tau>0,
\]
reveal the temporal dependence, when not being equal to $\E\left((x_t-\E(x_t))^2\right)\E\left((x_{t+\tau}-\E(x_{t+\tau}))^2\right)$.
A general idea in ICA is to find  for the standardized series $\bo x_t^{st}$ the orthogonal transformation which maximizes a selected objective function. In this paper we propose using a weighted sum of the squared linear and quadratic autocovariances of the components as the objective function.

\section{ARMA-GARCH model}\label{sec:ARMAGARCH}

Let $\mathcal{F}_t$ denote the information of a univariate time series process $(x_t)_{t=0,\pm 1,\pm 2,\dots}$ at time $t$. The conditional mean process of a time series process
\[
x_t=\E(x_t|\mathcal{F}_{t-1})+ z_t
\]
is often modeled as an ARMA$(p,q)$ process
\[
\E(x_t|\mathcal{F}_{t-1})=\sum_{i=1}^p\phi_i x_{t-i}+\sum_{j=1}^q\theta_j z_{t-j},
\]
where $\phi_1,\dots,\phi_p$ are the autoregressive coefficients and $\theta_1,\dots,\theta_q$ are the moving average coefficients. The ARMA process is causal if it can be written as a MA$(\infty)$ process
\[
x_t =\sum_{i=0}^{\infty}\psi_i z_{t-i}.
\]

The process $(z_t)_{t=0,\pm 1,\pm 2,\dots}$ is a white noise process and often it is assumed to consist of iid Gaussian random variables. The ARMA-GARCH process (see for example \citet{Weiss1984, Bollerslev1986}) is obtained when $z$ is a GARCH$(P,Q)$ process
\[
z_t = \sigma_t \epsilon_t,
\]
where $\epsilon_t$ is a Gaussian white noise process with unit variance and $\sigma_t^2$ is a conditional variance process
\begin{align*}
\sigma_t^2 &= \var(z_t| \mathcal{F}_{t-1} )= \omega+\sum_{i=1}^P\alpha_i z_{t-i}^2 +\sum_{j=1}^Q \beta_j \sigma^2_{t-j},
\end{align*}
with $\omega>0$ and $\alpha_i, \beta_j \ge 0$ for all $i=1,\dots,P$ and $j=1,\dots,Q$. The process is second-order stationary if $\sum_{i=1}^p \alpha_i + \sum_{j=1}^q \beta_j < 1$ \citep{Bollerslev1986}.

Our focus in this paper is on $p$-variate causal ARMA-GARCH$(1,1)$ processes $\bo x=(\bo x_t)_{t=0,\pm 1,\pm 2,\dots}$ with mutually independent components
\begin{align}
\label{armagarch}
x_{tj} &=\sum_{i=0}^{\infty}\psi_{ij} z_{t-i,j},\ \ j=1,\dots,p,
\end{align}
where $\sum_{i=0}^{\infty}\psi_{ij}^2=1$ and $z_{tj}$ is a GARCH$(1,1)$ process
\begin{align*}
z_{tj} &= \sigma_{tj} \epsilon_{tj},\ \ j=1,\dots,p,
\end{align*}
where $\epsilon_{tj}$ is a Gaussian white noise process with unit variance and
\begin{align*}
\sigma_{tj}^2 &= \omega_j +\alpha_j z_{t-1,j}^2 + \beta_j \sigma^2_{t-1,j},
\end{align*}
with $\alpha_j, \beta_j \ge 0$ and $\omega_j=1-\alpha_j-\beta_j$, which implies that $\var(z_{tj})=1$, and then also $\var(x_{tj})=1$. We restrict to processes with finite eighth moments. Notice that a GARCH$(1,1)$ process with parameters $\alpha$ and $\beta$ has finite eighth moments if and only if $\beta^4+4\beta^3\alpha+18\beta^2\alpha^2+60\beta\alpha^3+105\alpha^4<1$. For assessing the finiteness of even moments of a GARCH(1,1) process, see \citet{Bollerslev1986}.

\section{Generalized SOBI estimator}
\label{sec:GSOBI}

Let now $\bo S_\tau(\bo x)=\E(\bo x_t\bo x_{t+\tau}^T)$ denote the autocovariance matrix with lag $\tau$ for a centered process $\bo x$. Then, for a set of lags $\mathcal{T}=\{\tau_1,\ldots,\tau_K\}$, the rotation matrix $\bo U$ in the SOBI estimator~\citep{Belouchranietal1997} was defined as the orthogonal matrix which makes the transformed autocovariance matrices
\[
\bo U\bo S_{\tau_1}(\bo x_t^{st})\bo U^T,\dots,\bo U\bo S_{\tau_K}(\bo x_t^{st})\bo U^T
\]
as diagonal as possible, the diagonality being measured using the sum of squares of the off-diagonal elements, or equivalently the sum of squares of the diagonal elements since the sum of squares of all elements of $\bo U\bo S_{\tau_1}(\bo x_t^{st})\bo U^T$ is the same for any orthogonal $\bo U$. Hence, the SOBI unmixing matrix estimator $\bo \Gamma=(\bo\gamma_1,\dots,\bo\gamma_p)^T$ can be alternatively defined as the maximizer of the objective function
\[
\sum_{j=1}^p\sum_{k=1}^K\left(\E\left(\bo \gamma_j^T \bo x_t\bo \gamma_j^T\bo x_{t+\tau_k}\right)\right)^2,
\]
under the constraint $\bo\Gamma\bo S(\bo x)\bo\Gamma^T=\bo I_p$. Note that when $K=1$ the estimator is known as the AMUSE (algorithm for multiple unknown signals extraction) estimator~\citep{Tongetal1990} and can be obtained as the solution of an eigenvector-eigenvalue problem. Using more lags simultaneously, as in SOBI, is however preferred as AMUSE depends heavily on the choice of the lag. For general statistical properties of AMUSE and SOBI (especially assuming linear processes), computational issues, choice of lags and variants see
\citet{Miettinenetal2012,Miettinenetal2014,Miettinenetal2016,ILLNER2015,TASKINEN2016}.

A more flexible use of temporal dependence in ICA has been considered in \cite{Hyvarinen2001}, \cite{Shietal2009}, \cite{MattesonTsay2011} and \cite{Matilainenetal2017}, the last of which defined vSOBI estimator as the maximizer of
\[
\sum_{j=1}^p\sum_{k=1}^K \left(\E\left(G(\bo\gamma_j^T\bo x_t)G(\bo\gamma_j^T\bo x_{t+\tau_k})\right) - \E\left(G(\bo\gamma_j^T\bo x_t)\right)^2\right)^2,
\]
under the constraint $\bo\Gamma\bo S(\bo x)\bo\Gamma^T=\bo I_p$, where $G$ is a twice continuously differentiable function. The SOBI estimator is obtained by the choice $G(x)=x$.~\citet{Matilainenetal2017} also discussed the possibility of combining SOBI and vSOBI with another function $G$, and in this paper we study this combination method when $G(x)=x^2$ is chosen. We call the method as the generalized SOBI (gSOBI), and it maximizes
\begin{align}
\label{ObjFunc}
&b \sum_{\tau \in \mathcal{T}_1} \sum_{j=1}^p\left(\E\left(\bo \gamma_j^T \bo x_t\bo \gamma_j^T\bo x_{t+\tau})\right)\right)^2 \nonumber \\
&+(1-b)\sum_{\tau \in \mathcal{T}_2} \sum_{j=1}^p\left(\E\left(\left(\bo\gamma_j^T \bo x_t\right)^2\left(\bo\gamma_j^T\bo x_{t+\tau}\right)^2-1\right)\right)^2,\ \ 0\leq b\leq 1,
\end{align}
under the constraint $\bo\Gamma\bo S(\bo x)\bo\Gamma^T=\bo I_p$. Here $\mathcal{T}_1$ and $\mathcal{T}_2$ are the sets of lags for the linear and quadratic parts, respectively, and $b$ is a parameter which gives the weight of the linear part.

\subsection{Estimating equations}

The estimating equations for SOBI and vSOBI were derived in~\citet{Miettinenetal2016} and~\citet{Matilainenetal2017}, respectively. In a similar fashion, the estimating equations for gSOBI can be derived and they are as follows.
\begin{definition}
	\label{estimatingequations}
	The generalized SOBI unmixing matrix functional $\bo \Gamma=(\bo\gamma_1,\dots,\bo\gamma_p)^T$ solves the estimating equations
	\begin{align*}
	&\bo\gamma_j^T\bo T(\bo\gamma_l)=\bo\gamma_l^T\bo T(\bo\gamma_j) \ \ \text{ and } \ \ \bo\gamma_j^T\bo S\bo\gamma_l=\delta_{jl},\ \ j,l=1,\dots,p,
	\end{align*}
	where $\delta_{jl}=1$, if $j=l$ and 0 otherwise,
	\[
	\bo T(\bo\gamma)=b\bo T^s(\bo \gamma)+(1-b)\bo T^v(\bo \gamma)
	\]
	and
	\begin{align*}
	\bo T^s(\bo \gamma)=& \sum_{\tau \in \mathcal{T}_1}\left(\E\left(\bo\gamma^T\bo x_t\bo\gamma^T\bo x_{t+\tau}\right)\E\left(\left(\bo\gamma^T\bo x_{t+\tau}\right)\bo x_t+ \left(\bo\gamma^T\bo x_t\right)\bo x_{t+\tau}\right)\right), \\
	\bo T^v(\bo \gamma)=& \, 2\sum_{\tau \in \mathcal{T}_2}\left(\E\left(\left(\bo\gamma^T\bo x_t\right)^2\left(\bo\gamma^T\bo x_{t+\tau}\right)^2-1\right)\times\right. \\
	&\left. \E\left(\left(\bo\gamma^T\bo x_t\right)\left(\bo\gamma^T\bo x_{t+\tau}\right)^2\bo x_t+\left(\bo\gamma^T\bo x_t\right)^2\left(\bo\gamma^T\bo x_{t+\tau}\right)\bo x_{t+\tau}\right)\right).
	\end{align*}
\end{definition}

The estimating equations allow us to study the asymptotical properties of the gSOBI estimator.

\subsection{Asymptotic properties}

Due to affine equivariance of the gSOBI estimator, see~\cite{Miettinenetal2016} and~\cite{Matilainenetal2017}, the limiting distribution of $\sqrt{n}\,(\hat{\bo\Gamma}\bo\Omega-\bo I_p)$ is the same for any mixing matrix $\bo\Omega$. Hence, we may first assume that $\bo\Omega=\bo I_p$. The limiting distribution in general case then follows as  $\mbox{vec}(\hat{\bo\Gamma}\bo\Omega)=(\bo\Omega^T\otimes \bo I_p)\mbox{vec}(\hat{\bo\Gamma})$. Here vec is an operator which vectorizes an matrix by stacking its columns on top of each other.

Let now $\bo e_j$ denote a $p$-vector whose $j$th element equals to one and the others are zeros. Then we write
\begin{align*}
\bo T_j^s=&\sum_{\tau \in \mathcal{T}_1}\left(\sum_{t=1}^{n-\tau}\left(\bo e_j^T\bo x_t\bo e_j^T\bo x_{t+\tau}\right) \sum_{t=1}^{n-\tau}\left(\left(\bo e_j^T\bo x_t\right)\bo x_{t+\tau}+\left(\bo e_j^T\bo x_{t+\tau}\right)\bo x_t\right)\right), \\
\bo T_j^v=& \,2\sum_{\tau \in \mathcal{T}_2}\left(\sum_{t=1}^{n-\tau}\left((\bo e_j^T\bo x_t)^2(\bo e_j^T\bo x_{t+\tau})^2-1\right)) \times \right.  \\
&\left. \sum_{t=1}^{n-t}\left(\left(\bo e_j^T\bo x_t\right)\left(\bo e_j^T\bo x_{t+\tau}\right)^2\bo x_t+\left(\bo e_j^T\bo x_t\right)^2\left(\bo e_j^T\bo x_{t+\tau}\right)\bo x_{t+\tau}\right)\right), \\
\bo T_j=& \,b\bo T_j^s+(1-b)\bo T_j^v.
\end{align*}

For our main theorem, we make three assumptions on the independent components $\bo z$:
\begin{itemize}
	\item[(A1)] $\E(\bo z)=\bo 0$ and $\cov(\bo z)=\bo I_p$.
	\item[(A2)] The $p$ time series in $\bo z$ are independent.
	\item[(A3)] For each pair $j\neq l$ of independent components 
	\begin{itemize}
		\item[(a)] there is a $\tau \in \mathcal{T}_1$ such that $\mu_{\tau j}\neq \mu_{\tau l}$, or
		\item[(b)] there is a $\tau \in \mathcal{T}_2$ such that $\nu_{\tau j}^2+\nu_{\tau l}^2-2(\nu_{\tau j}+\nu_{\tau l})\mu_{\tau j}\mu_{\tau_l}\neq 0$,
	\end{itemize}
    where $\mu_{\tau j}=\E(z_{tj}z_{t+\tau,j})$ and $\nu_{\tau j}=\E(x_{tj}^2x_{t+\tau,j}^2)-1$.
\end{itemize}

The assumption (A3) is sufficient as such when $b<1$, and part (a) is required for $b=1$. In ARMA model, it can be shown that $\nu_{\tau j}=2\mu_{\tau j}^2$, and therefore, if $\mu_{\tau j}=\mu_{\tau l}$ for all $\tau \in \mathcal{T}_2$, part (b) of (A3) is not satisfied. Hence, identical ARMA processes without stochastic volatility are not allowed. On the other hand, identical (ARMA-)GARCH processes can be separated as well as mutually dissimilar processes. For $b=1$, Theorem~\ref{theorem1} and Corollary~\ref{cor1} generalize from ARMA~\citep{Miettinenetal2016} to ARMA-GARCH model. For $b\neq 1$ those results have not been derived before in any model.

The following results are proved in the Supplementary Material.
\begin{theorem}
	\label{theorem1}
	Under the assumptions (A1)-(A3) and $\bo\Omega=\bo I_p$, we have that $\hat{\bo\Gamma}\to_P \bo I_p$ and
	\begin{align*}
	&\sqrt{n}\,(\hat\gamma_{jj}-1)=-\frac{1}{2}\sqrt{n}(\hat{\bo S}_{jj}-1)+o_p(1), \\
	&\sqrt{n}\,\hat{\gamma}_{jl}=\left(\bo e_l^T\sqrt{n}\,\bo T_j-\bo e_j^T\sqrt{n}\,\bo T_l+\left(2b\sum_{\tau \in \mathcal{T}_1}\mu_{\tau j}(\mu_{\tau l}-\mu_{\tau j}) \right. \right. \\
	&\left. \left. +4(1-b)\sum_{\tau \in \mathcal{T}_2}\left(\nu_{\tau l}-\nu_{\tau j}(\nu_{\tau j}+1)+2\nu_{\tau l}\mu_{\tau j}\mu_{\tau l}\right)\right)\sqrt{n}\,\hat{\bo S}_{jl}\right) \\
	&\times \left(2b\sum_{\tau \in \mathcal{T}_1}\left(\mu_{\tau j}-\mu_{\tau l}\right)^2+4(1-b)\sum_{\tau \in \mathcal{T}_2}\left(\nu_{\tau j}^2+\nu_{\tau l}^2-2(\nu_{\tau j}+\nu_{\tau l})\mu_{\tau j}\mu_{\tau_l}\right)\right)^{-1} \\
	&+o_p(1), \ \ \ j\neq l,
	\end{align*}
	where $\hat{\bo S}$ is the sample covariance matrix.
\end{theorem}
From Theorem~\ref{theorem1} we obtain a straightforward corollary.
\begin{corollary}
	\label{cor1}
	In addition to (A1)-(A3), assume that $\bo z$ has finite eighth moments, and that the elements of $\hat{\bo S}$, $\bo T_j^s$, and $\bo T_j^v$, $j=1,\dots,p$ are asymptotically normally distributed. Then the limiting distribution of $\mbox{vec}(\hat{\bo\Gamma}-\bo I_p)$ is normal with zero mean and the following variances.
	\begin{align*}
	&ASV(\hat{\gamma}_{jj})=4^{-1}\left(\E\left[ x_{tj}^4\right]+\sum_{k=-\infty}^\infty\left(\E\left[x_{tj}^2x_{t+k,j}^2\right]-1\right)-1\right) \\
	&ASV\left(\hat{\gamma}_{jl}\right) = \left(b^2ASV\left(\bo e_l^T \bo T_j^s\right)+2b(1-b)  ASCOV\left(\bo e_l^T \bo T_j^s,\bo e_l^T \bo T_j^v\right) \right. \\
	&+(1-b)^2ASV\left(\bo e_l^T \bo T_j^v\right) +b^2ASV\left(\bo e_j^T \bo T_l^s\right)+2b(1-b)ASCOV\left(\bo e_j^T \bo T_l^s,\bo e_j^T \bo T_l^v\right) \\
	&+(1-b)^2ASV\left(\bo e_j^T \bo T_l^v\right)-2b^2ASCOV\left(\bo e_l^T \bo T_j^s,\bo e_j^T \bo T_l^s\right)-2(1-b)^2 \\
	&\times ASCOV\left(\bo e_l^T \bo T_j^v,\bo e_j^T \bo T_l^v\right)-2b(1-b)\left(ASCOV\left(\bo e_l^T \bo T_j^s,\bo e_j^T \bo T_l^v\right) \right. \\
	&\left.+ASCOV\left(\bo e_j^T \bo T_l^s,\bo e_l^T \bo T_j^v\right)\right) +2c_{jl}\left(b\,ASCOV\left(\bo e_l^T \bo T_j^s,\hat{\bo S}_{jl}\right)
	+(1-b) \right. \\
	&\times ASCOV\left(\bo e_l^T \bo T_j^v,\hat{\bo S}_{jl}\right)-b\,ASCOV\left(\bo e_j^T\bo T_l^s,\hat{\bo S}_{jl}\right)-(1-b) \\
	&\times \left. \left. ASCOV\left(\bo e_j^T\bo T_l^v,\hat{\bo S}_{jl}\right)\right) +c_{jl}^2ASV\left(\hat{\bo S}_{jl}\right)\right)\times\left(2b\sum_{\tau \in \mathcal{T}_1}\left(\mu_{\tau j}-\mu_{\tau l}\right)^2 \right. \\
	&\left.+4(1-b)\sum_{\tau \in \mathcal{T}_2}\left(\nu_{\tau j}^2+\nu_{\tau l}^2 -2(\nu_{\tau j}+\nu_{\tau l})\mu_{\tau j}\mu_{\tau l}\right)\right)^{-2}, \ \ \ j\neq l,
	\end{align*}
	where
	\[ c_{jl}=4(1-b)\sum_{\tau \in \mathcal{T}_2}\left(\nu_{\tau l}-\nu_{\tau j}(\nu_{\tau j}+1)+2\nu_{\tau l}\mu_{\tau j}\mu_{\tau l}\right)+2b\sum_{\tau \in \mathcal{T}_1}\left(\mu_{\tau j}(\mu_{\tau l}-\mu_{\tau j})\right).
	\]
\end{corollary}

If the components of $\bo z$ are ARMA-GARCH processes with finite eighth moments, the assumptions of Corollary~\ref{cor1} are satisfied. In the Supplementary material, we provide expression for the variances in the case of ARMA-GARCH$(1,1)$ processes.

\subsection{Tests for linear and quadratic autocorrelations}
\label{s:tests}

\citet{HuTsay2014} argued that especially in the context of econometric time series those latent series which exhibit ``more'' stochastic volatility are the most interesting ones. An ordering of components according to the ``degree'' of volatility clustering in the components is therefore an important task and will be considered next. We also demonstrate how to identify those components which do not contain such volatility clustering.

We suggest the following procedure for ordering the estimated components. Notice first that the straightforward use of
\begin{align}\label{ordering}
\sum_{\tau\in \mathcal{T}_2}\left(\sum_{t=1}^{n-\tau}\left(\hat{s}_{tj}^2\hat{s}_{t+\tau,j}^2\right)/(n-\tau)-1\right)^2,
\end{align}
as the criterion does not work, since linear autocorrelations imply quadratic autocorrelations. Therefore, at first, an ARMA model should be fitted to each time series exhibiting linear autocorrelation, after which the residual series can be plugged into~(\ref{ordering}) to order the components.

To test whether fitting the ARMA model is needed, we propose a modification of the Ljung-Box test, where the asymptotic variance of the test statistic is derived when assuming only symmetry and finite fourth moments of the time series.

\begin{definition}
	\label{LTest}
	The modified Ljung-Box test statistic using a set of lags $\mathcal{T}$ is given by
	\[
	L=n\sum_{\tau \in \mathcal{T}}\left(\sum_{t=1}^{n-\tau} \frac{x_tx_{t+\tau}}{n-\tau}\right)^2/V_{\tau},
	\]
	where $x$ is standardized to have zero mean and unit variance, and
	\[
	V_\tau=\sum_{t=1}^{n-\tau}\frac{x_t^2x_{t+\tau}^2}{n-\tau}+2\sum_{k=1}^{n-\tau-1}\frac{n-k}{n}\sum_{t=1}^{n-k-\tau}\frac{x_tx_{t+\tau}x_{t+k}x_{t+k+\tau}}{n-k-\tau}.
	\]
\end{definition}
Notice that, when $x_1,\dots,x_n$ are iid, $V_\tau\to_P 1$ and the classical Ljung-Box test is obtained by fixing $V_{\tau}=1$. However, in the case that $x$ has volatility clustering, $V_\tau>1$, and the estimation of $V_\tau$ is required to achieve the correct size of the test.

\begin{theorem}
	Under the null hypothesis that $\E(x_tx_{t+\tau})=0$ for all $\tau=1,2,\dots$, the limiting distribution of $L$ is the chi-squared distribution with $|\mathcal{T}|$ degrees of freedom.
\end{theorem}
Since $\E(x_tx_{t+\tau}x_{t+k}x_{t+k+\tau})$ is typically very close to zero for large $k$, we take the sum in the formula of $V_\tau$ only up to 20 in our simulations.

To test whether the quadratic autocorrelations are typical for homoscedastic time series, we use the following test statistic

\begin{definition}
	\label{QTest}
	The test statistic for volatility clustering, using a set of lags $\mathcal{T}$, is
	\[
	Q=n\sum_{\tau \in \mathcal{T}}\left(\sum_{t=1}^{n-\tau} \frac{x_t^2x_{t+\tau}^2}{n-\tau}-1\right)^2/4,
	\]
	where $x$ is standardized to have zero mean and unit variance.
\end{definition}
Before using the test statistic $Q$, one should remove the linear autocorrelation from the series. This can be performed, for example, by fitting an ARMA model and then applying the test statistic to the residual series.
\begin{theorem}
	Under the null hypothesis that $x_1,\dots,x_n$ are iid Gaussian, the distribution of $Q$ is the chi-squared distribution with $|\mathcal{T}|$ degrees of freedom.
\end{theorem}

The validity of the tests with finite sample sizes is studied in Section~\ref{sec:simulations}.

\section{Principal volatility component analysis and ICA}
\label{sec:pvca}

The lag-$l$ generalized kurtosis matrix of a $p$-variate time series $\bo x$ was defined in~\cite{HuTsay2014} as
\[
\bo g_l(\bo x)=\sum_{i=1}^p\sum_{j=i}^p\cov^2(\bo x_t\bo x_t^T,y_{t-l,ij}),
\]
where $y_{t-l,ij}$ is a function of $x_{t-l,i}x_{t-l,j}$ for $1\leq i,j\leq p$ and $l\geq 0$. For what follows, we restrict ourselves to the identity function $y_{t-l,ij}=x_{t-l,i}x_{t-l,j}$, and run the index $j$ from 1 to $p$, i.e., we define
\[
\bo g_l(\bo x)=\sum_{i=1}^p\sum_{j=1}^p\cov^2(\bo x_t\bo x_t^T,x_{t-l,i}x_{t-l,j}).
\]

Note that this definition is slightly different than the one originally suggested by \cite{HuTsay2014} where the second index starts only at $i$. Equivariance as proven below is however only achievable with this modification.

The cumulative generalized kurtosis matrix is then
\[
\bo G_m(\bo x)=\sum_{l=1}^m \bo g_l(\bo x).
\]
The following theorem states that this version of $\bo G_m$ possesses two properties that are required when used for ICA, see \cite{Ojaetal2006}. The proofs are given in the Supplementary Material.
\begin{theorem}
	\label{pvcthm1}
	(i) $\bo G_m$ is orthogonal equivariant ,i.e., for any $p$-variate random vector $\bo x$ and any $p\times p$ orthogonal matrix $\bo V$,
	\[
	\bo G_m(\bo V\bo x)=\bo V\bo G_m(\bo X)\bo V^T.
	\]
	(ii) $\bo G_m$ has independence property \citep{NordhausenTyler:2015}, i.e., if $\bo s$ has independent components, then $\bo G_m(\bo s)$ is diagonal.
\end{theorem}
Hence, we may construct the PVC estimator for the independent component problem as a simultaneous diagonalizer of the covariance matrix $\bo S$ and the cumulative generalized kurtosis matrix $\bo G_m$.
\begin{definition}
	The PVC unmixing matrix functional $\Gamma$ satisfies
	\[
	\bo\Gamma\bo S(\bo x)\bo\Gamma^T=\bo I_p \text{ and } \bo\Gamma\bo G_m(\bo x)\bo\Gamma^T=\bo D,
	\]
	where $\bo D$ is a diagonal matrix.
\end{definition}

This unmixing matrix functional is affine equivariant as one of the two matrices is affine equivariant and the other one is orthogonal equivariant~\citep{Ojaetal2006}.  
In practice, the estimator can be obtained by finding the eigendecomposition of $\bo G_m(\bo x^{st})$, which is the strategy of \citet{HuTsay2014} who recommended first whitening 
the times series using the covariance matrix and then computing the kurtosis matrix $\bo G_m$. Simultaneous diagonalization separates a pair of independent components only if the corresponding eigenvalues are distinct. In the case of PVC, we have to assume that for $i \neq j$
\[
\sum_{l=1}^m\E^2((s_{ti}^2-1)(s_{t-l,i}^2-1)) \neq \sum_{l=1}^m\E^2((s_{tj}^2-1)(s_{t-l,j}^2-1)).
\]
When the eigenvalues are distinct, simultaneous diagonalization of two matrices is a computationally simple method to solve the independent component problem, but the estimation efficiency is known to be inferior to more advanced methods.

\section{Simulation study}
\label{sec:simulations}

All following computations have been made in R \citep{R} using the packages JADE \citep{JADE}, tsBSS \citep{tsBSS}, forecast \citep{forecast} and fGarch \citep{fGarch}. In this simulation study gSOBI is compared to PVC. The vSOBI was compared to some other methods in~\citet{Matilainenetal2017}.

\subsection{Performance index}

To measure the performance of the estimates in the simulation study, we use the minimum distance index (MDI)~(\cite{Ilmonenetal2010})
\begin{align*}
\hat{D}= D(\hat{\bo\Gamma}\bo\Omega)=\frac{1}{\sqrt{p-1}}\inf_{\bo C\in
	\mathcal{C}}\|\bo C\hat{\bo\Gamma}\bo\Omega-\bo I_p\|,
\end{align*}
where $\|\cdot\|$ is the matrix (Frobenius) norm and
\[
\mathcal C=\{\bo C\ :\ \mbox{each row and column of $\bo C$ has exactly one non-zero
	element}  \}.
\]
The index takes values between zero and one, and zero indicates perfect separation. The minimum distance index is not difficult to compute and it has some attractive properties. It is invariant with respect to the change of the mixing matrix, and when the estimate $\hat{\bo\Gamma}$ is asymptotically normally distributed with expected value $\bo I_p$, then the limiting expected value of $n(p-1)\hat{D}^2$
is the sum of the limiting variances of the off-diagonal elements of $\sqrt{n}\,\hat{\bo\Gamma}$.

\subsection{Models}

First, in model (i), we consider a three-variate ARMA(1,1)-GARCH(1,1) model with the parameter values given in Table~\ref{table1}. Model (ii) is three-variate pure ARMA(1,1) and model (iii) is three-variate pure GARCH(1,1), where $\phi$ and $\theta$ of model (ii) and $\alpha$, $\beta$ and $\omega$ of model (iii) take values as in Table~\ref{table1}. In model (iv) there are two pure ARMA(1,1) and two pure GARCH(1,1) components with parameters as in $s_1$ and $s_2$ in Table~\ref{table1}. Since all estimators are affine equivariant, we use $\bo \Omega = \bo I_p$.

\begin{table}[htb]
	\centering
	\caption{The parameter values of the ARMA-GARCH time series in model (i). Parameter $\phi$ is the AR coefficient and $\theta$ is the MA coefficient.}
	\label{table1}
	\begin{tabular}{cccccc} \hline
		& $\alpha$ & $\beta$ & $\omega$ & $\phi$ & $\theta$  \\ \hline
		$s_1$ &  0.15    &  0.7    &  0.15    &  0.5   &  -0.1    \\
		$s_2$ &  0.1     &  0.8    &  0.1     &  0.2   &  0.8     \\
		$s_3$ &  0.05     &  0.9    &  0.05     &  0.1   &  0.1     \\ \hline
	\end{tabular}
\end{table}

\subsection{Efficiency results}

\begin{figure}
	\begin{center}
		\includegraphics[width=3in]{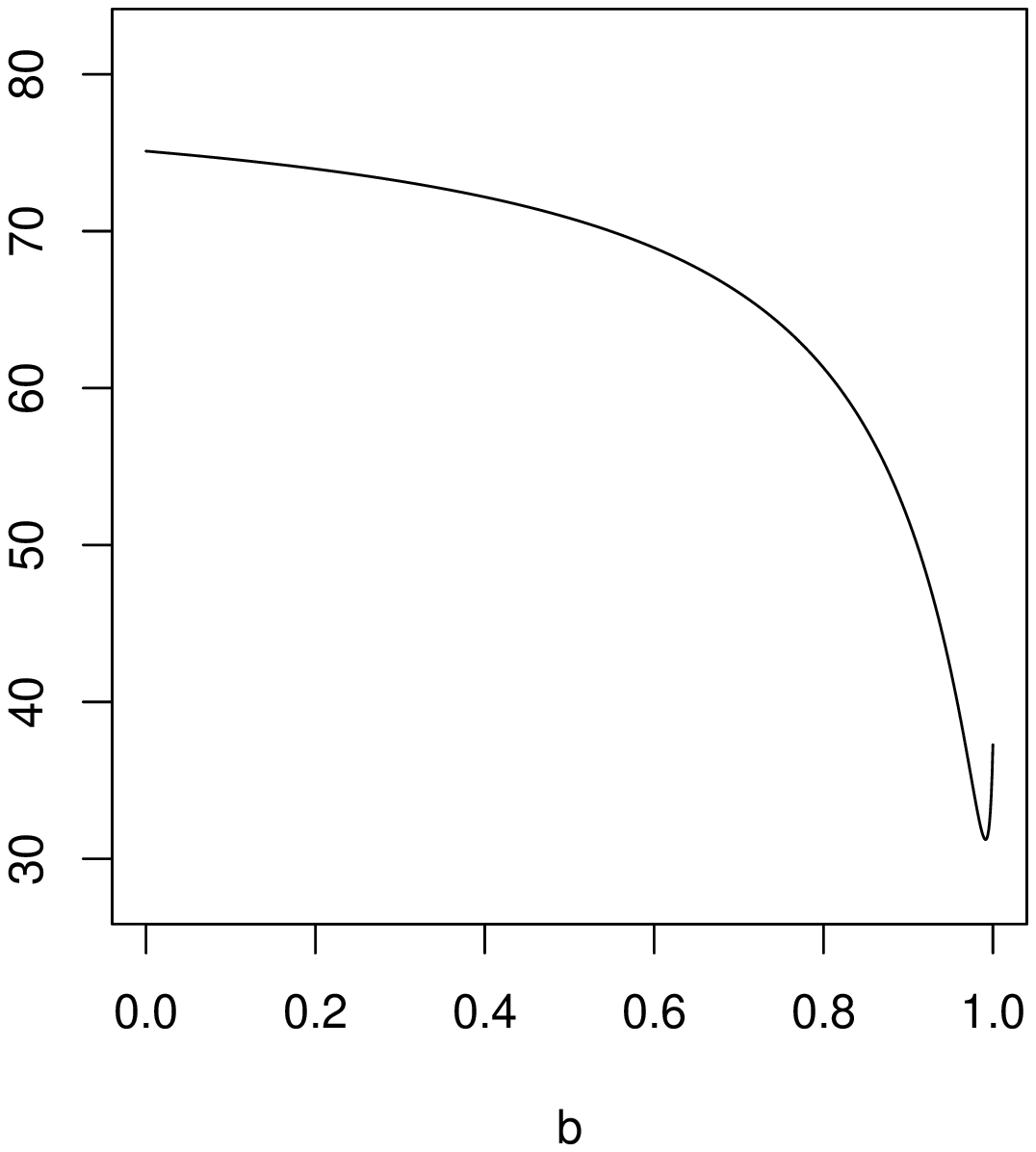}
		\includegraphics[width=3in]{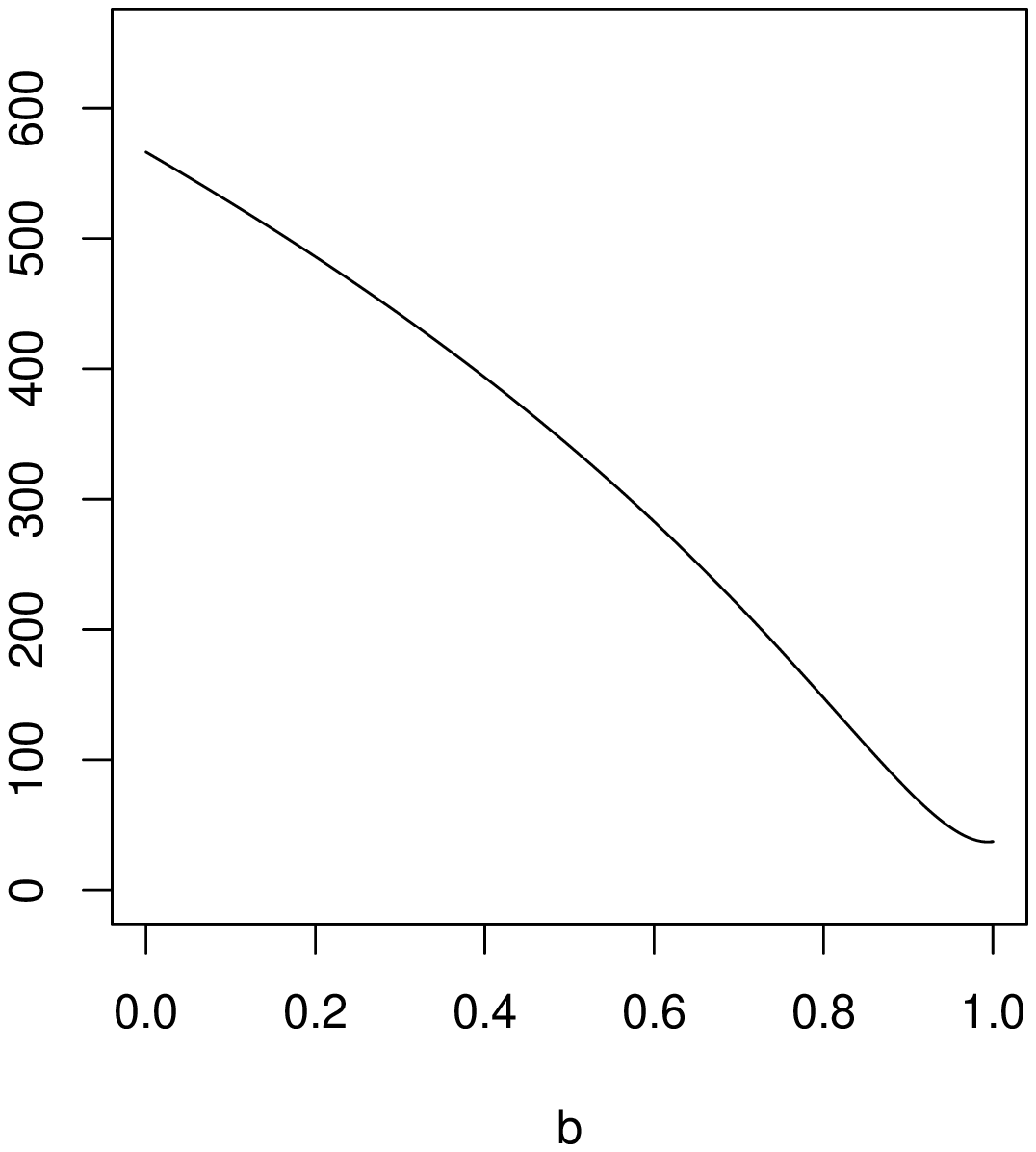} \\
		\includegraphics[width=3in]{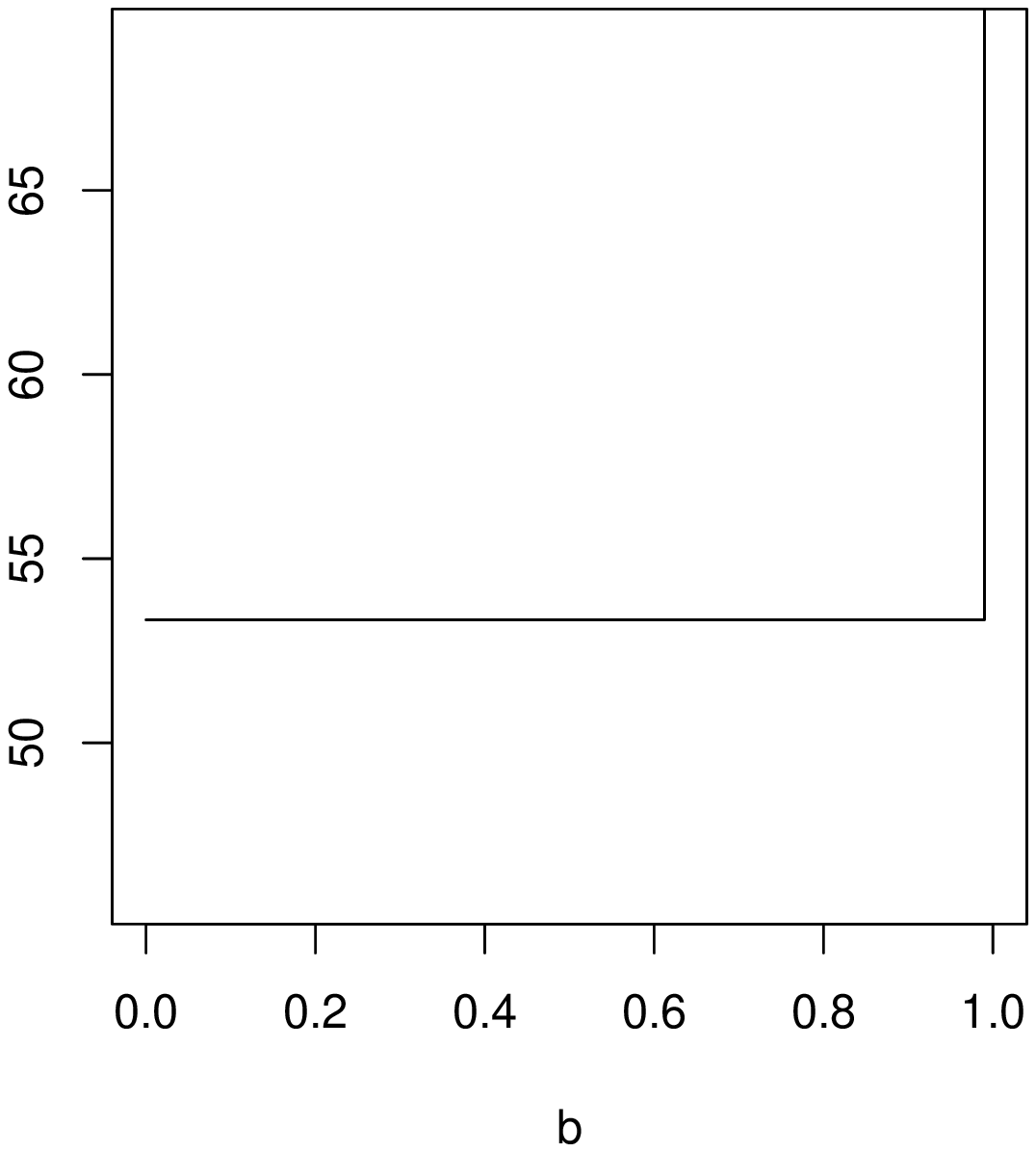}
		\includegraphics[width=3in]{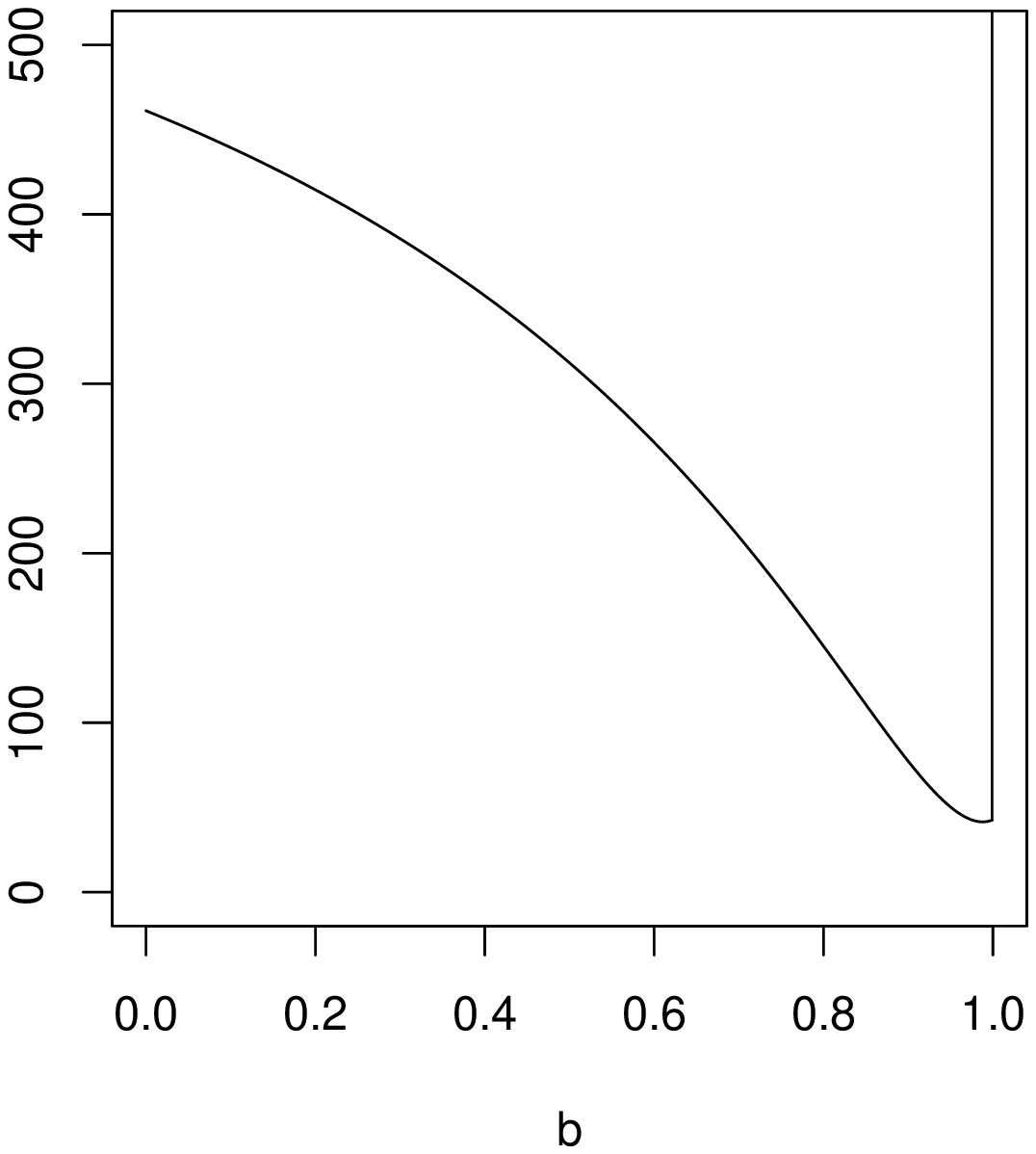}
	\end{center}
	\caption{The limiting expected value of $n(p-1)\hat{D}^2$ as function of $b$ in model (i) on the top left, (ii) on the top right, (v) on the bottom left and (iv) on the bottom right.\label{fig1}}
\end{figure}

In our simulations, we selected the lag sets $\mathcal{T}_1=\left\{1,2,3\right\}$ and $\mathcal{T}_2=\left\{1,2,3\right\}$, but used also lag set $\tilde{\mathcal{T}}_1=\left\{1,\dots,12\right\}$ for $b=0.9$. The limiting expected value of $n(p-1)\hat{D}^2$ as a function of $b$ in models (i)-(iv) is plotted in Figure~\ref{fig1}. In model (iii), the expected value is constant for $b<1$ and goes to infinity when $b=1$. Large values of $b$ are preferable in this asymptotic comparison, excluding SOBI $(b=1)$ which totally fails in models (iii) and (iv) where there are more than one components with zero linear autocorrelation for all lags.

\begin{figure}
	\begin{center}
		\includegraphics[width=2.7in]{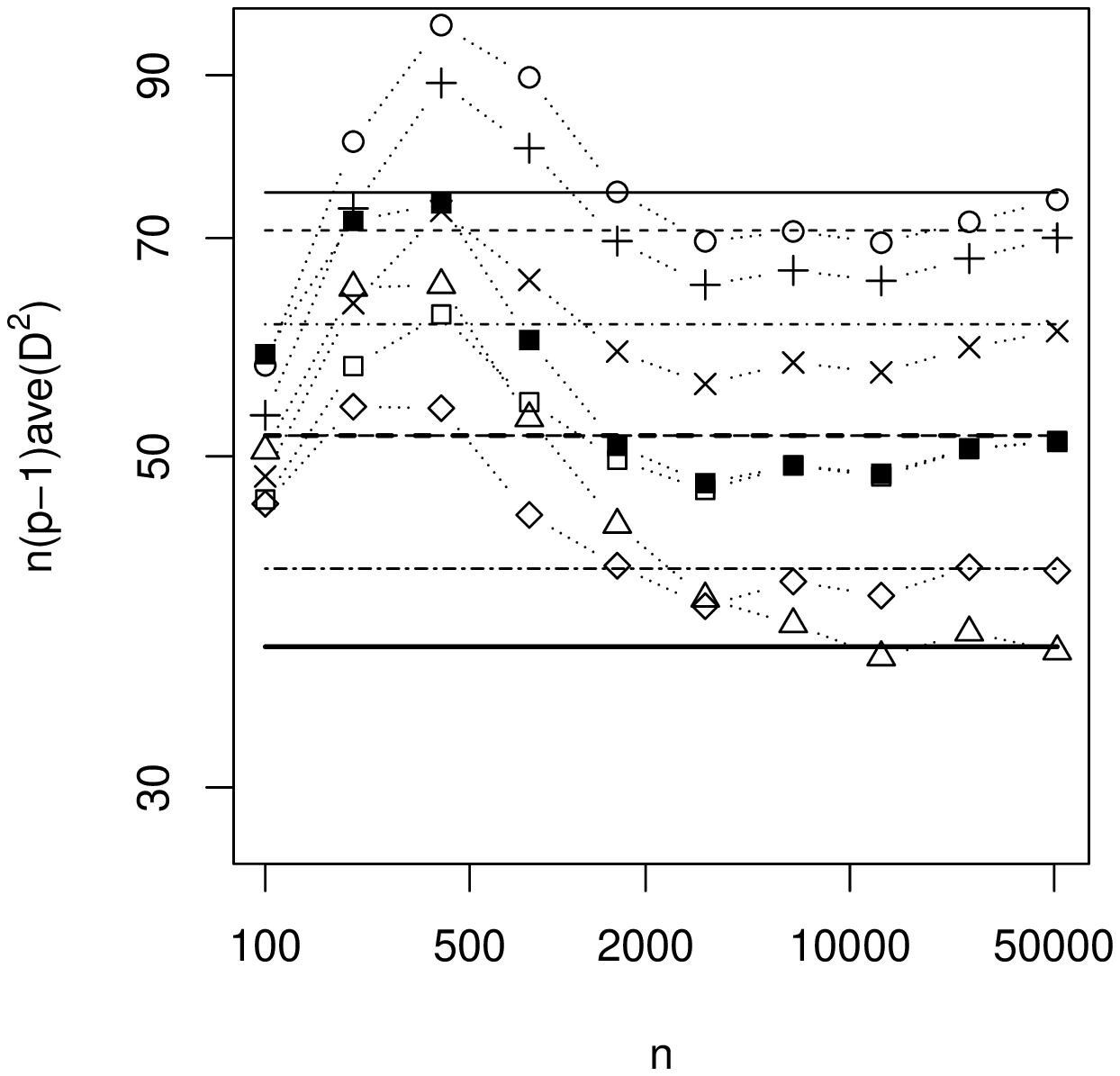}
		\includegraphics[width=3.3in]{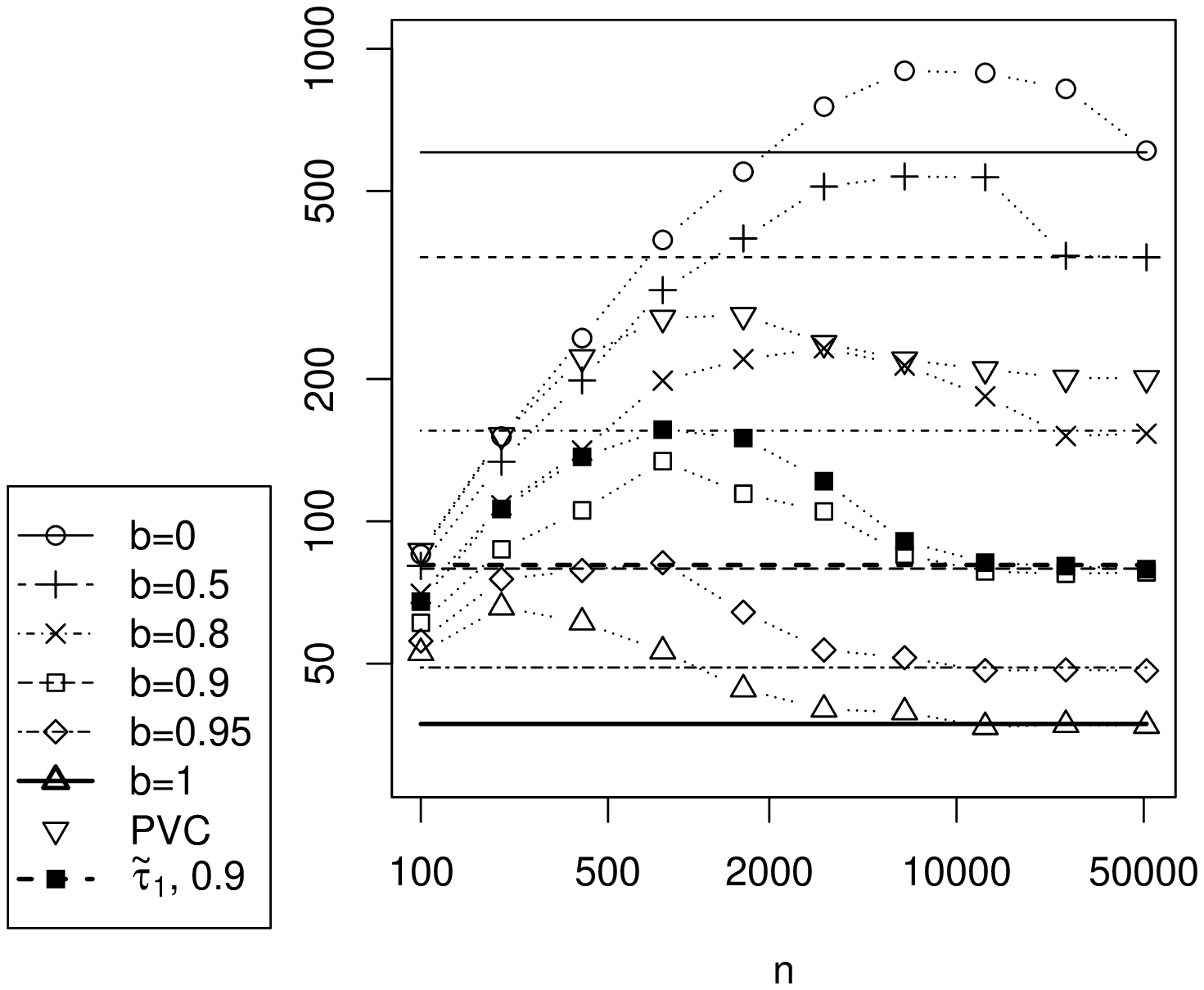} \\
		\includegraphics[width=2.7in]{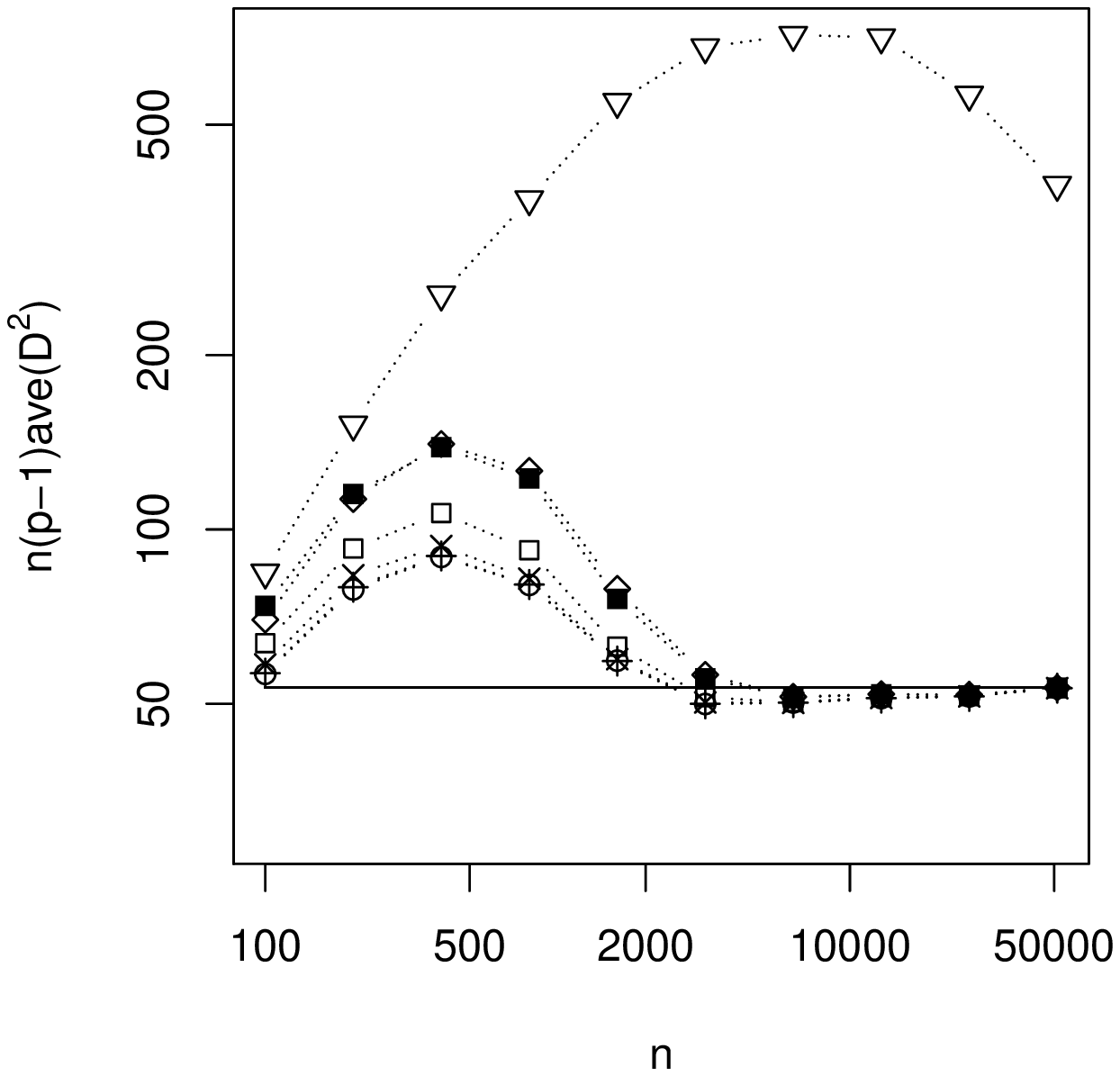}
		\includegraphics[width=3.3in]{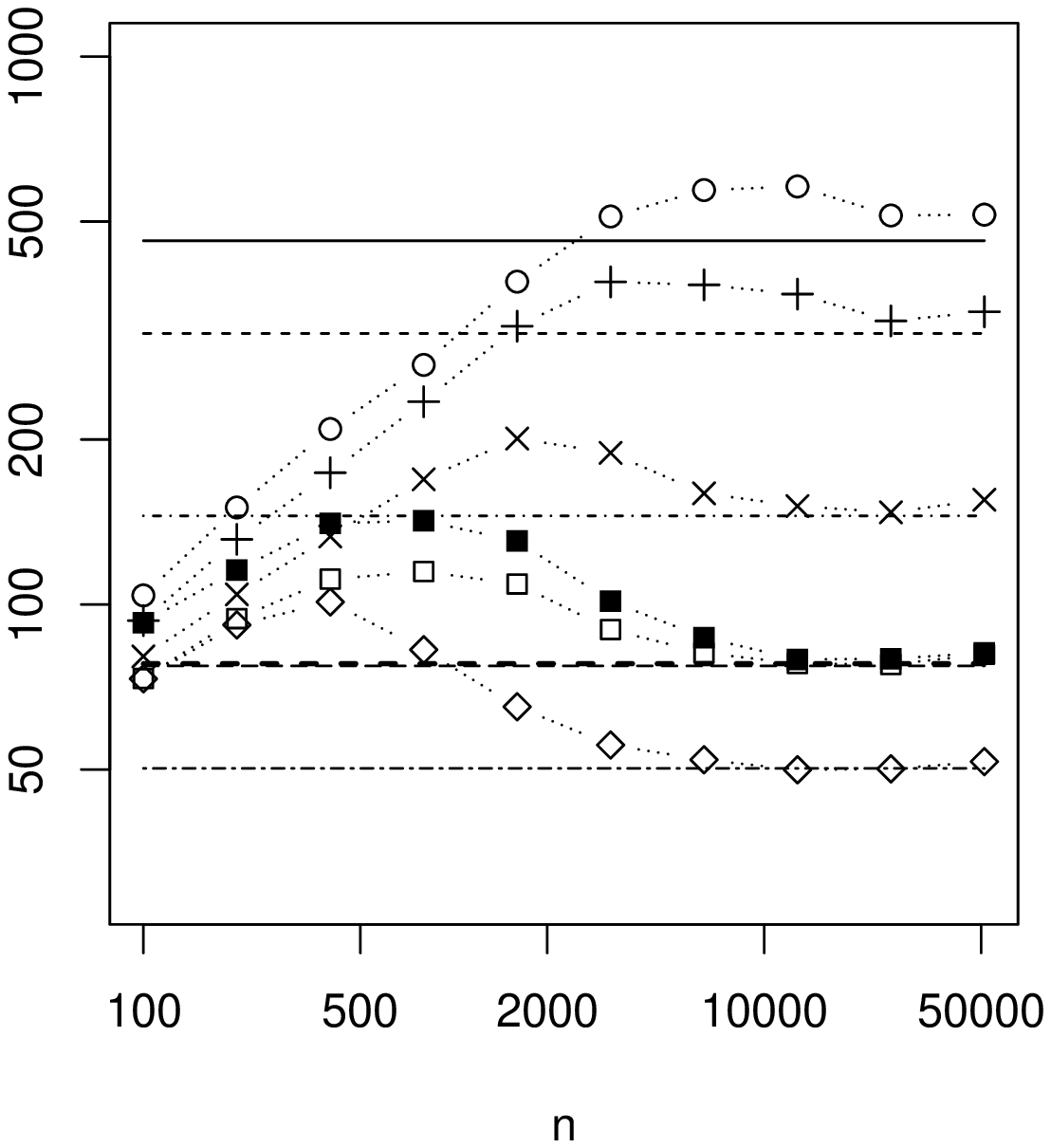}
	\end{center}
	\caption{Simulation results of models (i)--(iv). The points show the averages over 2000 repetitions for each sample size, and the horizontal lines give the asymptotic expected values. PVC and SOBI (b=1) are  out of the plotting region in models (iii) and (iv), and therefore not visible. Model (i) on the top left, (ii) on the top right, (iii) on the bottom left and (iv) on the bottom right. Both axes are on log scales.\label{fig2}}
\end{figure}

The finite sample behavior of the six gSOBI estimators given by $b=0$, $b=0.5$, $b=0.8$, $b=0.9$, $b=0.95$, and $b=1$, and the principal volatility component (PVC) estimator with $m=5$ is compared in a simulation study. For sample sizes $n=100,200,400,\dots,51200$ we generated 2000 realizations from each model. Figure~\ref{fig2} shows the averages of $n(p-1)\hat{D}^2$ over the 2000 repetitions. As expected, SOBI is very poor in models (iii) and (iv), and it is omitted from those figures. PVC is dropped out of the figures corresponding to models (i) and (iv) as the averages of $n(p-1)\hat{D}^2$ grow with $n$. In model (iii), all gSOBI estimators with $b<1$ are asymptotically equally efficient, and in the simulations they are practically equally good from $n=3200$ on. With smaller sample sizes the unnecessary use of linear autocorrelations cause loss in efficiency. However, $b=0.5$ is already as good as $b=0$. In model (ii), with only ARMA components, the quadratic autocorrelations are still nonzero and thus also $b=0$ works, but choices with $b$ closer to 1 are much better. All in all, it seems clear that the combinations are better than extreme choices $b=0$ and $b=1$, and that $b$ should be larger than 0.5. We thus suggest the use of value $b=0.9$. The lag set $\mathcal{T}_1$ is here slightly better than $\tilde{\mathcal{T}}_1$, but the latter is generally a safer choice. 

\subsection{Tests}

For each data set which were generated from models (i), (ii) and (iii), tests of linear and quadratic autocorrelations were performed as described in Section~\ref{s:tests}. In the tests, we used the estimates of the independent component time series given by gSOBI with $b=0.9$. Table~\ref{table2} presents the proportions of rejected null hypotheses at significance level 0.05, when testing the linear autocorrelations in models (i), (ii) and (iii), and Table~\ref{table3} presents the proportions of rejected null hypotheses in the quadratic autocorrelations test. In the test for linear autocorrelations, both the modified and the classical Ljung-Box tests were applied. The results in Table~\ref{table2} show that in the GARCH model (iii), where the linear autocorrelations are zero, the modified Ljung-Box test has the correct rate of false rejections, whereas the classical test rejects too often. On the other hand, the classical test seems to be more efficient.
Recall that the test for quadratic autocorrelations is applied to residual series of ARMA models, and that the estimated components are ordered according to the values of the test statistic. Therefore, especially in the case of model (ii), we are mainly interested in the sum of the rejection probabilities, which should be $0.15$ when there is no volatility clustering in the components. Table~\ref{table3} shows that the rate of false rejections is close to the desired level from $n=800$ upwards.

\begin{table}
	\caption{The proportion of rejections in the test of linear autocorrelation at significance level 0.05 in models (i), (ii) and (iii). The modified Ljung-Box on the left and the classical Ljung-Box on the right. \label{table2}}
	\begin{center}
	 \begin{tabular}{ccccc} \hline
	  Model & $n$ & $s_1$ & $s_2$ & $s_3$ \\ \hline
			& 100 & 0.768 \ $\slash$ \ 0.839 &  0.874 \ $\slash$ \ 0.898 & 0.486 \ $\slash$ \ 0.523 \\
			& 200 & 0.927 \ $\slash$ \ 0.945 &  0.948 \ $\slash$ \ 0.962 & 0.622 \ $\slash$ \ 0.665 \\
	    (i) & 400 & 0.988 \ $\slash$ \ 0.994 &  0.988 \ $\slash$ \ 0.994 & 0.893 \ $\slash$ \ 0.925 \\
			& 800 & 1 \ $\slash$ \ 1         &  0.999 \ $\slash$ \ 1     & 0.995 \ $\slash$ \ 0.998 \\
			& $\geq$ 1600 & 1 \ $\slash$ \ 1 &  1 \ $\slash$ \ 1         & 1 \ $\slash$ \ 1 \\ \hline
			& 100 & 0.833 \ $\slash$ \ 0.848 &  0.847 \ $\slash$ \ 0.867 & 0.528 \ $\slash$ \ 0.545 \\
			& 200 & 0.921 \ $\slash$ \ 0.920 &  0.919 \ $\slash$ \ 0.926 & 0.690 \ $\slash$ \ 0.696 \\
	   (ii) & 400 & 0.983 \ $\slash$ \ 0.986 &  0.980 \ $\slash$ \ 0.981 & 0.933 \ $\slash$ \ 0.934 \\
			& 800 & 1 \ $\slash$ \ 1         &  1 \ $\slash$ \ 1         & 0.997 \ $\slash$ \ 0.997 \\
			& $\geq$ 1600 & 1 \ $\slash$ \ 1 &  1 \ $\slash$ \ 1         & 1 \ $\slash$ \ 1 \\  \hline
			& 100 & 0.122 \ $\slash$ \ 0.203 &  0.127 \ $\slash$ \ 0.184 & 0.139 \ $\slash$ \ 0.178 \\
			& 200 & 0.094 \ $\slash$ \ 0.181 &  0.111 \ $\slash$ \ 0.175 & 0.097 \ $\slash$ \ 0.139 \\
			& 400 & 0.068 \ $\slash$ \ 0.158 &  0.089 \ $\slash$ \ 0.142 & 0.096 \ $\slash$ \ 0.130 \\
			& 800 & 0.059 \ $\slash$ \ 0.156 &  0.064 \ $\slash$ \ 0.134 & 0.065 \ $\slash$ \ 0.094 \\
	  (iii) & 1600 & 0.053 \ $\slash$ \ 0.147 & 0.063 \ $\slash$ \ 0.121 & 0.049 \ $\slash$ \ 0.080 \\
			& 3200 & 0.052 \ $\slash$ \ 0.148 & 0.056 \ $\slash$ \ 0.118 & 0.054 \ $\slash$ \ 0.086 \\
			& 6400 & 0.055 \ $\slash$ \ 0.153 & 0.048 \ $\slash$ \ 0.108 & 0.060 \ $\slash$ \ 0.087 \\
			& 12800 & 0.056 \ $\slash$ \ 0.153 & 0.054 \ $\slash$ \ 0.120 & 0.055 \ $\slash$ \ 0.092 \\
			& 25600 & 0.048 \ $\slash$ \ 0.150 & 0.053 \ $\slash$ \ 0.119 & 0.052 \ $\slash$ \ 0.082 \\
			& 51200 & 0.056 \ $\slash$ \ 0.144 & 0.046 \ $\slash$ \ 0.106 & 0.049 \ $\slash$ \ 0.072 \\ \hline
		\end{tabular}
	\end{center}
\end{table}

\begin{table}
	\caption{The proportion of rejections in the test of quadratic autocorrelations at 0.05 significance level in models (i), (ii) and (iii). \label{table3}}
	\begin{center}
		\begin{tabular}{cccccccccc} \hline
			Model &  & (i) &  &  & (ii) &  &  & (iii) & \\
			$n$   & $s_1$ & $s_2$ & $s_3$ & $s_1$ & $s_2$ & $s_3$ & $s_1$ & $s_2$ & $s_3$ \\ \hline
			100   & 0.292 & 0.180 & 0.155 & 0.073 & 0.075 & 0.077 & 0.443 & 0.308 & 0.242 \\
			200   & 0.526 & 0.322 & 0.258 & 0.067 & 0.063 & 0.062 & 0.669 & 0.500 & 0.336 \\
			400   & 0.841 & 0.596 & 0.361 & 0.062 & 0.064 & 0.059 & 0.878 & 0.715 & 0.440\\
			800   & 0.983 & 0.911 & 0.598 & 0.058 & 0.049 & 0.055 & 0.987 & 0.929 & 0.658 \\
			1600  & 1     & 0.996 & 0.880 & 0.056 & 0.049 & 0.058 & 1     & 1     & 0.880 \\
			3200  & 1     & 1     & 0.996 & 0.044 & 0.044 & 0.055 & 1     & 1     & 0.996 \\
			6400  & 1     & 1     & 1     & 0.053 & 0.054 & 0.048 & 1     & 1     & 1 \\
			12800 & 1     & 1     & 1     & 0.050 & 0.053 & 0.051 & 1     & 1     & 1 \\
			25600 & 1     & 1     & 1     & 0.054 & 0.050 & 0.049 & 1     & 1     & 1 \\
			51200 & 1     & 1     & 1     & 0.041 & 0.052 & 0.051 & 1     & 1     & 1 \\\hline
		\end{tabular}
	\end{center}
\end{table}

In models (i) and (iii), we also examined ordering the components according to volatility,
\[
\sum_{\tau\in \mathcal{T}_2}\left(\sum_{t=1}^{n-\tau}\left(r_{tj}^2r_{t+\tau,j}^2\right)/(n-\tau)-1\right)^2,
\]
where $r$ denotes the residual time series of the independent components after fitting the ARMA model. The expected values of these for the three components were 1.88, 1.69 and 0.53. Table~\ref{table4} shows how often the components are in the expected order at different sample sizes. In the GARCH model (iii), the order is quite consistently the expected one already at small values of $n$, whereas in the ARMA-GARCH model (i), larger sample size is needed for that.

\begin{table}
	\caption{The proportions of correct ordering of the components and the proportions of correct identification of the component with least volatility clustering, in models (i) and (iii).  \label{table4}}
	\begin{center}
		\begin{tabular}{ccccccccccc} \hline
			$n$ & 100  & 200  & 400  & 800  & 1600 & 3200 & 6400 & 12800 & 25600 & 51200 \\ \hline
		    (i) & 0.20 & 0.29 & 0.40 & 0.50 & 0.66 & 0.79 & 0.89 & 0.96  & 0.99  & 1.00  \\
		  (iii) & 0.72 & 0.284 & 0.95 & 0.99 & 1.00 & 1.00 & 1.00 & 1.00  & 1.00  & 1.00 \\ \hline
		\end{tabular}
	\end{center}
\end{table}

\section{Application}
\label{sec:data}

In this section we apply our proposed gSOBI method to weekly log returns of seven different exchange rates against the U.S. Dollar. We use the same data as in \cite{HuTsay2014} and compare our results to those derived from their principal volatility component (PVC) analysis.

The seven exchange rates in the original data, obtained from the database of \cite{WeeklyReturns}, are those of Australian Dollar (AUD), Canadian Dollar (CAD), Singaporean Dollar (SGD), British Pound (GBP), Swiss Franc (CHF), Norwegian Kroner (NOK) and Swedish Kroner (SEK). The data included are from March 11, 2000 to November 26, 2011.

Let $\bo y_t$ be a vector of average values of the exchange rates against U.S. Dollar for week $t$. As observed values for gSOBI algorithm we use the logarithmic returns of weekly averages, i.e. values $\bo x_t = \mbox{log}(\bo y_t) - \mbox{log}(\bo y_{t-1})$, $t = 1, 2, \ldots, T$. The total length of the series is $T = 605$.

First, we extract the mutually independent components $\bo s_t$ using gSOBI algorithm with $b = 0.9$, as recommended earlier, and $\mathcal{T}_1 = \{1,\dots,12\}$ and $\mathcal{T}_2 = \{1, 2, 3\}$. We fit stationary and non-seasonal ARMA($p,q$) models to each series where the test statistic $L$ indicates the presence of the linear autocorrelation. Here $p$ and $q$ are determined using AIC. Then the residuals of the ARMA fit are used in further analysis, in order to better detect the volatility clustering after removing linear autocorrelations.

\begin{figure}
	\begin{center}
		\includegraphics[width=5in]{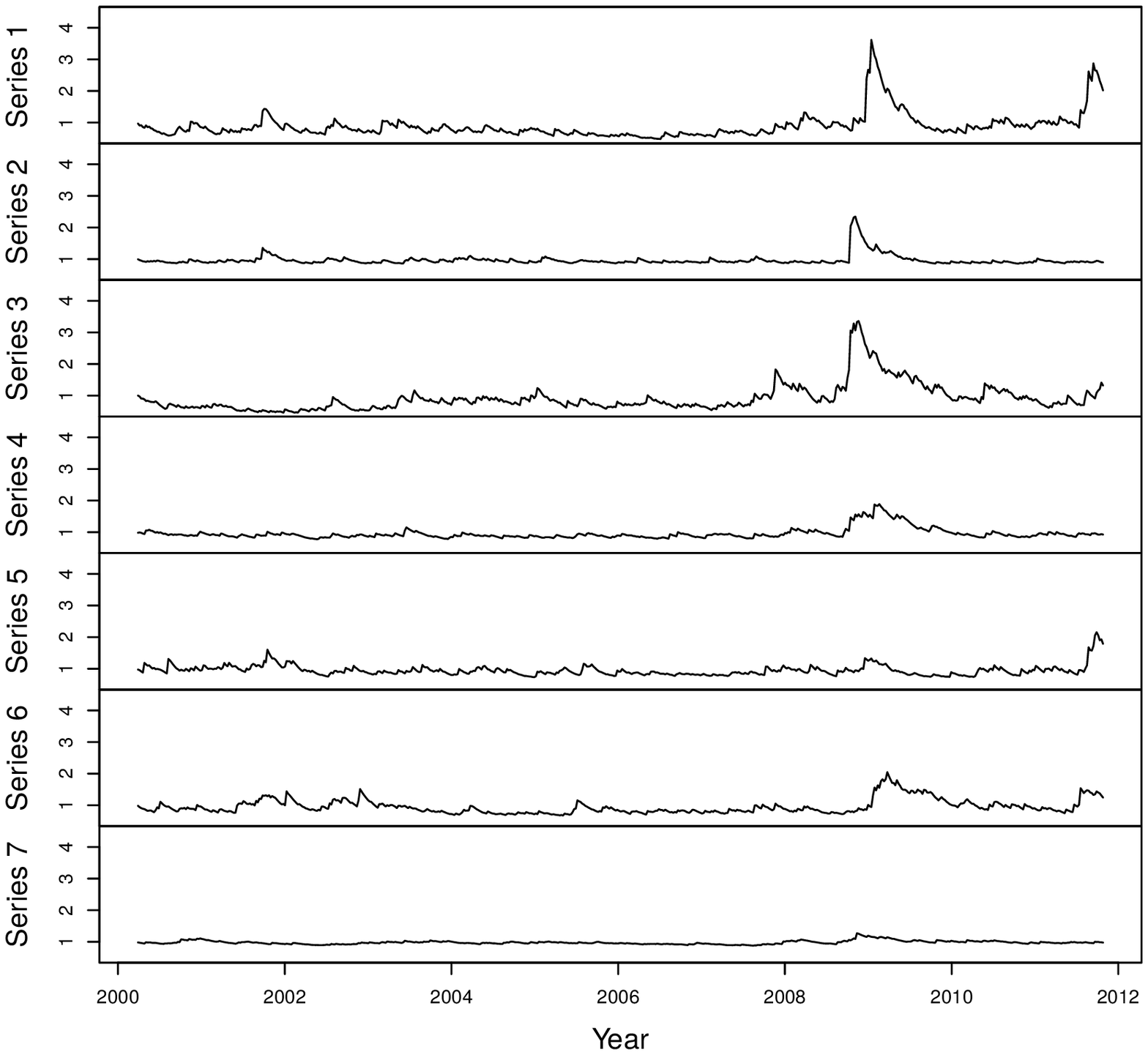}
	\end{center}
	\caption{Volatility components of $\bo s_t$ based on the exchange rate data.\label{fig::volatility}}
\end{figure}

For each residual series we test the hypothesis that the quadratic autocorrelations are zero using the test statistic $Q$. High value of a quadratic autocorrelation indicates that the volatility changes in time.
Components are then ordered according to the degree of volatility clustering, i.e., according to the value of $Q$.
For each series with stochastic volatility present we fit a univariate GARCH(1,1) model (as in \cite{HuTsay2014}) and calculate the volatility series based on that. Finally the volatility series are plotted in Figure \ref{fig::volatility}.

Table \ref{table::example} summarizes the results of the $L$- and $Q$-statistics and their significance for each independent component based on lags $1, \ldots, 5$. Also possible ARMA($p,q$) and GARCH(1,1) fits for each relevant component are included. Changing lags used in gSOBI or $L$- and $Q$-tests does not change the results significantly.

Most of the independent components have nonzero linear autocorrelation, and all of the components have statistically significant $Q$-statistic. Hence, a GARCH(1,1) is fitted to all components and the volatility series are plotted for each of the seven components. From Figure \ref{fig::volatility} it can be seen that series 7 has very little volatility clustering compared to all other series and other series have peaks of different magnitudes in volatility especially during the years 2008 and 2009.

\begin{table}\label{table::example}
	\caption{Values of $L$ and $Q$-statistic with their significance level and parameter estimates based on ARMA($p,q$) fit, when applicable, and GARCH(1,1) fits for each series. $^*$ = $Q$-test is performed on the residual series.}
	\begin{center}
		\begin{tabular}{ccccccccc}
			\hline
			&  Series 1 & Series 2 & Series 3 & Series 4 & Series 5 & Series 6 & Series 7 \\
			\hline
			$L$ & 14.33 & 8.09 & 6.01 & 15.19 & 21.14 & 13.01 & 23.79 \\
			$p$-value & 0.014 & 0.152 & 0.305 & 0.010 & $<$ 0.001 & 0.023 & $<$ 0.001\\
			ARMA($p,q$):  \\
			$\phi_1$ & -- & -- & -- & -- & -0.465 & 0.734 & -- \\
			$\theta_1$ & 0.271 & -- & -- & 0.187 & 0.668 & -0.610 & 0.199 \\
			$\theta_2$ & -- & -- & -- & -- & -- & -0.190 & -- \\
			\hline
			$Q$ & 7418.7$^*$ & 5242.1 & 4625.1 & 464.6$^*$ & 432.9$^*$ & 276.6$^*$ & 18.1$^*$ \\
			$p$-value & $<$ 0.001 & $<$ 0.001 & $<$ 0.001 & $<$ 0.001 & $<$ 0.001 & $<$ 0.001 & 0.0029\\
			GARCH(1,1): \\
			$\omega$ & 0.035 & 0.120 & 0.017 & 0.041 & 0.076 & 0.039 & 0.039 \\
			$\alpha$ & 0.136 & 0.051 & 0.141 & 0.050 & 0.087 & 0.086 & 0.021 \\
			$\beta$ & 0.849 & 0.820 & 0.847 & 0.905 & 0.833 & 0.875 & 0.938 \\ \hline
		\end{tabular}
	\end{center}
\end{table}

The key differences in the analyses of \cite{HuTsay2014} and ours are that they first fit a VAR model to remove linear autocorrelations and extract the PVCs from the residuals of the model, and that the components are ordered with different criteria, which are related to the estimation methods. Notice also that the eigenvectors in \cite{HuTsay2014} were not orthogonal, which means that the generalized kurtosis matrix was not symmetrized. This implies that the PVC's were not made uncorrelated, unlike in most of the ICA methods. In the results, we see that the series with most volatility clustering in our analysis, i.e. series 1 and series 2, have similar patterns to series $v1$ and $v3$ in \cite{HuTsay2014}, respectively. In the other end, both found one component with much less volatility changes than the other components. However, we still got statistically significant value of test statistic $Q$ for it, whereas \cite{HuTsay2014} found one component without volatility clustering.

\section{Conclusion}
\label{sec:conc}

Combined use of linear and quadratic autocorrelations in independent component analysis enables efficient estimation of latent time series of different types. The asymptotic distribution of the novel estimator was derived, and asymptotic variances were computed in four ARMA-GARCH models. The gSOBI estimator outperforms its two components, SOBI and vSOBI, especially when some, but not all, of the components have nonzero linear autocorrelations. Even if all the linear autocorrelations are zero, gSOBI with any weighting is asymptotically equally efficient as vSOBI.
Equal weights for the linear and quadratic parts are not recommended, because then the quadratic part dominates. We also derived new tests for checking whether linear or quadratic autocorrelations are zeros. The test statistic for linear autocorrelations is the same as in the well-known Ljung-Box test, but we consider the distribution of the test statistic under null hypothesis of linear independence instead of full independence.

\section*{Supplementary Material}

\begin{description}
	\item[Proofs:] Proofs of the theoretical results. (Appendix)
	
	\item[R-package tsBSS:] R-package tsBSS contains codes for the computation of the estimates and tests described in the article, and the dataset which was used in the article. (GNU zipped tar file)
	
\end{description}

\bibliographystyle{chicago}

\bibliography{refs3}

\section*{Appendix}
\subsection*{Proof of Theorem 1}
Write
\begin{align*}
&\hat{\bo T}(\hat{\bo \gamma})=b\sum_{\tau \in \mathcal{T}_1}\left(\frac{1}{n}\sum_{t=1}^{n-\tau}\left(\hat{\bo\gamma}^T\bo x_t\hat{\bo\gamma}^T\bo x_{t+\tau}\right) \frac{1}{n}\sum_{t=1}^{n-\tau}\left((\hat{\bo\gamma}^T\bo x_{t+\tau})\bo x_t+ (\hat{\bo\gamma}^T\bo x_t)\bo x_{t+\tau}\right)\right) \\
&+2(1-b)\sum_{\tau \in \mathcal{T}_2}\left(\frac{1}{n}\sum_{t=1}^{n-\tau}\left((\hat{\bo\gamma}^T\bo x_t)^2(\hat{\bo\gamma}^T\bo x_{t+\tau})^2-1\right)\times\right. \\
&\left. \frac{1}{n}\sum_{t=1}^{n-\tau}\left((\hat{\bo\gamma}^T\bo x_t)(\hat{\bo\gamma}^T\bo x_{t+\tau})^2\bo x_t+ (\hat{\bo\gamma}^T\bo x_t)^2(\hat{\bo\gamma}^T\bo x_{t+\tau})\bo x_{t+\tau}\right)\right) \ \ \text{ and}  \\
&\hat{\bo T}_j:=\hat{\bo T}(\hat{\bo\gamma}_j)
\end{align*}

The estimating equations of the gSOBI estimate are
\[
\hat{\bo\gamma}_l^T \hat{\bo T}_j=\hat{\bo\gamma}_j^T \hat{\bo T}_l\ \  \text{ and } \ \ \hat{\bo\gamma}_j^T S\hat{\bo\gamma}_l=\delta_{jl}.
\]

We begin with
\begin{align*}
\sqrt{n}\,\hat{\bo\gamma}_l^T \hat{\bo T}_j=&\sqrt{n}\,(\hat{\bo\gamma}_l-\bo e_l)^T\hat{\bo T}_j+ \sqrt{n}\,\bo e_l^T\left(\hat{\bo T}_j-\left(2b\sum_{\tau \in \mathcal{T}_1}\mu_{\tau j}^2 \right.\right. \\
&\left.\left. -4(1-b)\sum_{\tau \in \mathcal{T}_2}\nu_{\tau j}(\nu_{\tau j}+1)\right)\bo e_j\right)+o_P(1).
\end{align*}
Taylor's and Slutsky's theorems give
\begin{align*}
&\sqrt{n}\,\left(\hat{\bo T}_j-\left(2b\sum_{\tau \in \mathcal{T}_1}\mu_{\tau j}^2-4(1-w)\sum_{\tau \in \mathcal{T}_2}\nu_{\tau j}(\nu_{\tau j}+1)\right)\bo e_j\right) \\
=&\sqrt{n}\,\left(\bo T_j-\left(2b\sum_{\tau \in \mathcal{T}_1}\mu_{\tau j}^2-4(1-b)\sum_{\tau \in \mathcal{T}_2}\nu_{\tau j}(\nu_{\tau j}+1)\right)\bo e_j\right) \\
&+\left(2b\sum_{\tau \in \mathcal{T}_1}\left(\mu_{\tau j}^2\bo e_j\bo e_j^T+\mu_{\tau j}\E\left(\bo x_t\bo x_{t+\tau}+\bo x_{t+\tau}\bo x_t \right)\right)+16\sum_{\tau \in \mathcal{T}_2}\left((\nu_{kj}+1)^2\bo e_j\bo e_j^T \right. \right. \\
&+\nu_{\tau j}\E\left(2x_{t+k,j}^2\bo x_t\bo x_t^T+4x_{t,j}x_{t+k,j}\bo x_t\bo x_{t+k}^T +4x_{t,j}x_{t+\tau,j}\bo x_{t+\tau}\bo x_t^T \right. \\
&\left. \left. \left. +2x_{t,j}^2\bo x_{t+\tau}\bo x_{t+\tau}^T\right)\right)\right)\sqrt{n}\,(\hat{\bo\gamma}_j-\bo e_j)+o_P(1).
\end{align*}
After some calculations we obtain
\begin{align*}
\sqrt{n}\,\hat{\bo\gamma}_l^T \hat{\bo T}_j=& \left(2b\sum_{\tau \in \mathcal{T}_1}\mu_{\tau j}+4(1-b)\sum_{\tau \in \mathcal{T}_2}\nu_{kj}(\nu_{\tau j}+1)\right)\sqrt{n}\,\hat{\gamma}_{lj}+\bo e_l^T\sqrt{n}\,\bo T_j \\
&+\left(2b\sum_{\tau \in \mathcal{T}_1}\mu_{\tau j}\mu_{\tau l}+4\sum_{\tau \in \mathcal{T}_2}\nu_{\tau j}\left(1+2\mu_{\tau j}\mu_{\tau l} \right)\right)\sqrt{n}\,\hat{\gamma}_{jl}+o_P(1).
\end{align*}
Finally, the estimating equations yield the equation
\begin{align*}
&\left(2b\sum_{\tau \in \mathcal{T}_1}\mu_{\tau j}^2+4(1-b)\sum_{\tau \in \mathcal{T}_2}\nu_{\tau j}(\nu_{\tau j}+1)\right)\sqrt{n}\,\left(-\hat{\gamma}_{jl}-\hat{S}_{jl}\right)+\bo e_l^T\sqrt{n}\,\bo T_j \\
&+\left(2b\sum_{\tau \in \mathcal{T}_1}\mu_{\tau j}\mu_{\tau l}+4(1-b)\sum_{\tau \in \mathcal{T}_2}\nu_{\tau j}\left(1+2\mu_{\tau j}\mu_{\tau l}\right)\right)\sqrt{n}\,\hat{\gamma}_{jl}+o_P(1)\\
=& \left(2b\sum_{\tau \in \mathcal{T}_1}\mu_{\tau l}^2+4(1-b)\sum_{\tau \in \mathcal{T}_2}\nu_{\tau l}(\nu_{\tau l}+1)\right)\sqrt{n}\,\hat{\gamma}_{jl}+\bo e_j^T\sqrt{n}\,\bo T_l\\
&+\left(2b\sum_{\tau \in \mathcal{T}_1}\mu_{\tau j}\mu_{\tau l}+4(1-b)\sum_{\tau \in \mathcal{T}_2}\nu_{\tau l}\left(1+2\mu_{\tau j}\mu_{\tau l}\right)\right)\sqrt{n}\,(-\hat{\gamma}_{jl}-\hat{S}_{jl}),
\end{align*}
and hence
\begin{align*}
\sqrt{n}\,\hat{\gamma}_{jl}=&\left(\bo e_l^T\sqrt{n}\,\bo T_j-\bo e_j^T\sqrt{n}\,\bo T_l+\left(2b\sum_{\tau \in \mathcal{T}_1}\mu_{\tau j}(\mu_{\tau l}-\mu_{\tau j}) \right. \right. \\
&\left. \left. +4(1-b)\sum_{\tau \in \mathcal{T}_2}\left(\nu_{\tau l}-\nu_{\tau j}(\nu_{\tau j}+1)+2\nu_{\tau l}\mu_{\tau j}\mu_{\tau l}\right)\right)\sqrt{n}\,\hat{S}_{jl}\right) \\
&\times \left(2b\sum_{\tau \in \mathcal{T}_1}\left(\mu_{\tau j}-\mu_{\tau l}\right)^2+4(1-b)\sum_{\tau \in \mathcal{T}_2}\left(\nu_{\tau j}^2+\nu_{\tau l}^2-2(\nu_{\tau j}+\nu_{\tau l})\mu_{\tau j}\mu_{\tau_l}\right)\right)^{-1}.
\end{align*}

\subsection*{Asymptotic variances and covariances for Corollary 1}

Assuming an ARMA-GARCH$(1,1)$ model with finite eighth moments, we have the following asymptotic variances. The required expected values of the process $\bo x$ are given below.
\begin{align*}
&ASV(\hat{S}_{jl}) = 1+2\sum_{k=1}^K\left(\sum_{i=0}^{\infty}\psi_{ij}\psi_{i+k,j}\sum_{i=0}^{\infty}\psi_{il}\psi_{i+k,l}\right) \\
&ASV(\hat{S}_{jj}) = \E\left[ x_{tj}^4\right]+2\sum_{k=1}^\infty\E\left[x_{tj}^2x_{t+k,j}^2 \right]-2 \\
&ASV(\bo e_l^T\bo T_j^v) = 4\sum_{k\in\mathcal{T}_2}\sum_{i\in\mathcal{T}_2}\nu_{kj}\nu_{ij}\sum_{r=-\infty}^\infty\left(\E\left[x_{tj}x_{t+k,j}^2x_{t+r,j}x_{t+r+i,j}^2\right]\E\left[x_{tl}x_{t+r,l}\right] \right. \\
&+\E\left[x_{tj}^2x_{t+k,j}x_{t+r,j}x_{t+r+i,j}^2\right]\E\left[x_{t+k,l}x_{t+r,l}\right]+\E\left[x_{tj}x_{t+k,j}^2 x_{t+r,j}^2x_{t+r+i,j}\right] \\
&\left. \times\E\left[x_{tl}x_{t+r+i,l}\right]+\E\left[x_{tj}^2x_{t+k,j} x_{t+r,j}^2x_{t+r+i,j}\right]\E\left[x_{t+k,l}x_{t+r+i,l}\right]\right) \\
&ASV(\bo e_l^T\bo T_j^s) = \sum_{k\in\mathcal{T}_1}\sum_{i\in\mathcal{T}_1}\mu_{kj}\mu_{ij}\sum_{r=-\infty}^\infty\left(\mu_{r+i-k,j}\mu_{rl}+\mu_{r-k,j}\mu_{r+i,l} \right. \\
&\left. +\mu_{r+i,j}\mu_{r-k,l}+\mu_{rj}\mu_{r+i-k,l}\right) \\
&ASCOV\left(\bo e_l^T \bo T_j^v,\hat{S}_{jl}\right) = 2\sum_{k\in\mathcal{T}_2}\nu_{kj}\sum_{r=-\infty}^\infty\left(\E\left[x_{tj}x_{t+k,j}^2x_{t+r,j}\right]\mu_{rl}\right. \\
&\left. +\E\left[x_{tj}^2x_{t+k,j}x_{t+r,j}\right]\mu_{r-k,l}\right) \\
&ASCOV\left(\bo e_l^T \bo T_j^s,\hat{S}_{jl}\right) =\sum_{k\in\mathcal{T}_2}\mu_{kj}\sum_{r=-\infty}^\infty\left(\mu_{r-k,j}\mu_{rl}+\mu_{rj}\mu_{r-k,l}\right) \\
&ASCOV\left(\bo e_l^T\bo T_j^v,\bo e_j^T\bo T_l^v\right) = 4\sum_{k\in\mathcal{T}_2}\sum_{i\in\mathcal{T}_2}\nu_{kj}\nu_{il}\sum_{r=-\infty}^\infty\left(\E\left[x_{tj}x_{t+k,j}^2x_{t+r,j}\right] \right. \\
&\times \E\left[x_{t,l}x_{t+r,l}x_{t+r+i,l}^2\right]+\E\left[x_{tj}^2x_{t+k,j}x_{t+r,j}\right]\E\left[x_{t+k,l}x_{t+r,l}x_{t+r+i,l}^2\right] \\
& +\E\left[x_{tj}x_{t+k,j}^2x_{t+r+i,j}\right] \E\left[x_{tl}x_{t+r,l}^2x_{t+r+i,l}\right] \\
&\left. +\E\left[x_{tj}^2x_{t+k,j}x_{t+r+i,j}\right]\E\left[x_{t+k,l}x_{t+r,l}^2x_{t+r+i,l}\right] \right)
\end{align*}

\begin{align*}
&ASCOV\left(\bo e_l^T\bo T_j^s,\bo e_j^T\bo T_l^s\right) = \sum_{k\in\mathcal{T}_1}\sum_{i\in\mathcal{T}_1}\mu_{kj}\mu_{il}\sum_{r=-\infty}^\infty\left(\mu_{r-k,j}\mu_{r+i,l}+\mu_{r+i-k,j}\mu_{rl} \right. \\
&\left. +\mu_{rj}\mu_{r+i-k,l}+\mu_{r+i,j}\mu_{r-k,l}\right), \\
&ASCOV\left(\bo e_l^T\bo T_j^s,\bo e_l^T\bo T_j^v\right) = 2\sum_{k\in\mathcal{T}_1}\sum_{i\in\mathcal{T}_2}\mu_{kj}\nu_{ij}\sum_{r=-\infty}^\infty\left(\E\left[x_{tj}x_{t+r,j}x_{t+r+i,j}^2\right]\mu_{r-k,l} \right. \\
&+\E\left[x_{tj}x_{t+r,j}^2x_{t+r+i,j}\right]\mu_{r+i-k,l}+\E\left[x_{t+k,j}x_{t+r,j}x_{t+r+i,j}^2\right]\mu_{rl} \\
&\left.+\E\left[x_{t+k,j}x_{t+r,j}^2x_{t+r+i,j}\right]\mu_{r+i,l}\right), \\
&ASCOV\left(\bo e_j^T\bo T_l^s,\bo e_l^T\bo T_j^v\right) =2 \sum_{k\in\mathcal{T}_1}\sum_{i\in\mathcal{T}_2}\mu_{kl}\nu_{ij}\sum_{r=-\infty}^\infty\left(\E\left[x_{t+k,j}x_{t+r,j}x_{t+r+i,j}^2\right]\mu_{rl} \right. \\
&+\E\left[x_{t+k,j}x_{t+r,j}^2x_{t+r+i,j}\right]\mu_{r+i,l}+\E\left[x_{tj}x_{t+r,j}x_{t+r+i,j}^2\right]\mu_{r-k,l} \\
&\left.+\E\left[x_{tj}x_{t+r,j}^2x_{t+r+i,j}\right]\mu_{r+i-k,l}\right).
\end{align*}

\subsection*{Expected values of an ARMA-GARCH$(1,1)$ process}

Assume that $x$ and $z$ are from the ARMA-GARCH$(1,1)$ model~(1). Then $\E\left[x_t \right]=\E\left[z_t \right]=0$ and $\E\left[x_t^2 \right]=\E\left[z_t^2 \right]=1$, and
\begin{align*}
\E\left[z_t^4\right] =\, & \frac{3\left(\omega^2+2\omega\left(\alpha+\beta\right)\right)}{1-3\alpha^2-2\alpha\beta-\beta^2}, \\
\E\left[z_t^6\right] =\, & \frac{15\left(\omega^2\left(\omega+3\alpha+3\beta\right)+\omega\E\left[z_t^4\right]\left(3\alpha^2+2\alpha\beta+\beta^2\right)\right)}{1-15\alpha^3-9\alpha^2\beta-3\alpha\beta^2-\beta^3}, \\
\E\left[z_t^2z_{t+\tau}^2\right] =\, & \omega\sum_{k=0}^{\tau-1}\left(\alpha+\beta\right)^k+\left(\alpha+\beta\right)^{\tau-1}\left(\alpha+\beta/3\right)\E\left[z_t^4\right],\ \ \tau>0,  \\
\E\left[z_t^4z_{t+\tau}^2\right] =\, & \omega\sum_{k=0}^{\tau-1}\left(\alpha+\beta\right)^k\E\left[z_t^4\right]+\left(\alpha+\beta\right)^{\tau-1}\left(\alpha+\beta/5\right)\E\left[z_t^6\right],\ \ \tau>0, \\
\E\left[z_t^2z_{t+1}^4\right] =\, & 3\left(\omega^2+2\omega\left(\alpha+\beta/3\right)+\left(\alpha^2+2\alpha\beta/5+\beta^2/15\right)\E\left[z_t^6\right]\right), \\
\E\left[z_t^2z_{t+\tau}^4\right] =\, & \omega^2\sum_{k=0}^{\tau-1}3^{k+1}\left(\alpha^2+2\alpha\beta/3+\beta^2/3\right)^k+2\omega\left(\alpha+\beta\right)\sum_{k=0}^{\tau-2}3^{k+1} \times\\
&\left(\alpha^2+2\alpha\beta/3+\beta^2/3\right)^k\E\left[z_t^2z_{t+\tau-k-1}^2\right]+2\omega\left(\alpha+\beta/3\right)\E\left[z_t^4\right] \\
&+3^\tau\left(\alpha^2+2\alpha\beta/3+\beta^2/3\right)^{\tau-1}\left(\alpha^2+2\alpha\beta/5+\beta^2/15\right)\E\left[z_t^6\right], \\
&\tau>1,
\end{align*}

\begin{align*}
\E\left[z_t^2z_{t+\tau_1}^2z_{t+\tau_2}^2\right]=\, &\omega^2\left(\sum_{i=0}^{\tau_2-\tau_1-1}\sum_{j=0}^{\tau_1-1}c_{i+j,\alpha,\beta}+\sum_{i=1}^{\tau_1}\sum_{j=i}^{\tau_1}\left(1+\delta_{\{j>i\}} c_{\tau_2-\tau_1-1+j-i,\alpha,\beta}h_{i,\alpha,\beta} \right)\right) \\
&+\omega\E\left[z_t^4\right]\left(\alpha+\frac{\beta}{3}\right)\left(\sum_{i=1}^{\tau_1}2c_{\tau_2-1-i}h_{i,\alpha,\beta}+\sum_{i=1}^{\tau_2-\tau_1}2c_{\tau_1-2+i}\right)+ \\
&+\E\left[z_t^6\right]\left(\alpha^2+\frac{2\alpha\beta}{5}+\frac{\beta^2}{15}\right)c_{\tau_2-\tau_1-1,\alpha,\beta}h_{\tau_1,\alpha,\beta},
\end{align*}
where $\delta$ is the Kronecker delta and
\[
c_{\tau,\alpha,\beta}=\sum_{k=0}^\tau \alpha^k\beta^{\tau-k}\binom{\tau}{k}, \ \ \ h_{\tau,\alpha,\beta}=\sum_{i=0}^{\tau-1}\sum_{k=0}^{\tau}\sum_{r=0}^{\tau}3^r\binom{\tau-1}{i}\binom{i+1}{r}\binom{\tau-i-1}{k-r}.
\]

\begin{align*}
&\E\left[x_t^6\right]= \sum_{j=0}^\infty \psi_j^6\E\left[z_t^6\right]+\sum_{j=0}^\infty \sum_{k\neq j} 15\psi_j^4\psi_k^2\E\left[z_t^4z_{t+j-k}^2\right] \\
&+\sum_{j=0}^\infty\sum_{k\neq j}\sum_{l\neq j,l\neq k} 15 \psi_j^2\psi_k^2\psi_l^2\E\left[z_t^2z_{t+j-k}^2z_{t+j-l}^2\right] \\
&\E\left[x_t^4\right]=\sum_{j=0}^\infty \psi_j^4\E\left[z_t^4\right]+\sum_{j=0}^\infty \sum_{k\neq j} 3\psi_j^2\psi_k^2\E\left[z_t^2z_{t+j-k}^2\right] \\
&\E\left[x_t^4x_{t+\tau}^2\right]=\sum_{j=0}^\infty \psi_j^4\psi_{j+\tau}^2\E\left[z_t^6\right]+\sum_{j=0}^\infty \sum_{k\neq j+\tau} \psi_j^4\psi_k^2\E\left[z_t^4z_{t+j-k+\tau}^2\right] \\
&+\sum_{j=0}^\infty \sum_{k\neq j} \left(6\psi_j^2\psi_{j+\tau}^2\psi_k^2+8\psi_j^3\psi_{j+\tau}\psi_k\psi_{k+\tau}\right)\E\left[z_t^4z_{t+j-k}^2\right] \\
&+\sum_{j=0}^\infty\sum_{k\neq j}\sum_{l\neq j+\tau,l\neq k+\tau} 3\psi_j^2\psi_k^2\psi_l^2\E\left[z_t^2z_{t+j-k}^2z_{t+j-l+\tau}^2\right] \\
&+\sum_{j=0}^\infty\sum_{k\neq j}\sum_{l\neq j,l\neq k} 12\psi_j^2\psi_k\psi_{k+\tau}\psi_l\psi_{l+\tau}\E\left[z_t^2z_{t+j-k}^2z_{t+j-l}^2\right] 
\end{align*}

\begin{align*}
&\E\left[x_t^4x_{t+\tau_1}x_{t+\tau_2}\right]=\sum_{j=0}^\infty \psi_j^4\psi_{j+\tau_1}\psi_{j+\tau_2}\E\left[z_t^6\right]+\sum_{j=0}^\infty \sum_{k\neq j+\tau_1} \psi_j^4\psi_k\psi_{k+\tau_2-\tau_1}  \\
&\times \E\left[z_t^4z_{t+j-k+\tau_1}^2\right] +\sum_{j=0}^\infty \sum_{k\neq j}\left(4\psi_j^3\psi_{j+\tau_1}\psi_k\psi_{k+\tau_2}+4\psi_j^3\psi_{j+\tau_2}\psi_k\psi_{k+\tau_1} \right. \\
&\left.+6\psi_j^2\psi_{j+\tau_1}\psi_{j+\tau_2}\psi_k^2\right)\E\left[z_t^4z_{t+j-k}^2\right]+\sum_{j=0}^\infty\sum_{k\neq j}\sum_{l\neq j+\tau_1,l\neq k+\tau_1} 3\psi_j^2\psi_k^2\psi_l \psi_{l+\tau_2-\tau_1}  \\
&\times \E\left[z_t^2z_{t+j-k}^2z_{t+j-l+\tau_1}^2\right]+\sum_{j=0}^\infty\sum_{k\neq j}\sum_{l\neq j,l\neq k} 12\psi_j^2\psi_k\psi_{k+\tau_1}\psi_l\psi_{l+\tau_2}  \\
&\times \E\left[z_t^2z_{t+j-k}^2z_{t+j-l}^2\right] \\
&\E\left[x_t^3x_{t+\tau}^3\right]=\sum_{j=0}^\infty \psi_j^3\psi_{j+\tau}^3\E\left[z_t^6\right]+\sum_{j=0}^\infty \sum_{k\neq j+\tau} 3\psi_j^3\psi_{j+\tau}\psi_k^2\E\left[z_t^4z_{t+j-k+\tau}^2\right] \\
&+\sum_{j=0}^\infty \sum_{k\neq j}\left(9\psi_j^2\psi_{j+\tau}^2\psi_k\psi_{k+\tau}+3\psi_j\psi_{j+\tau}^3\psi_k^2\right)\E\left[z_t^4z_{t+j-k}^2\right] \\
&+\sum_{j=0}^\infty\sum_{k\neq j}\sum_{l\neq j+\tau,l\neq k+\tau} 9\psi_j^2\psi_k\psi_{k+\tau}\psi_l^2\E\left[z_t^2z_{t+j-k}^2z_{t+j-l+\tau}^2\right] \\
&+\sum_{j=0}^\infty\sum_{k\neq j}\sum_{l\neq j,l\neq k} 6\psi_j\psi_{j+\tau}\psi_k\psi_{k+\tau}\psi_l\psi_{l+\tau}\E\left[z_t^2z_{t+j-k}^2z_{t+j-l}^2\right] \\
&\E\left[x_t^3x_{t+\tau_1}^2x_{t+\tau_2}\right]=\sum_{j=0}^\infty \psi_j^3\psi_{j+\tau_1}^2\psi_{j+\tau_2}\E\left[z_t^6\right]+\sum_{j=0}^\infty \sum_{k\neq j+\tau_1} \left(\psi_j^3\psi_{j+\tau_2}\psi_k^2 \right. \\
&\left.+2\psi_j^3\psi_{j+\tau_1}\psi_k\psi_{k+\tau_2-\tau_1}\right)\E\left[z_t^4z_{t+j-k+\tau_1}^2\right]+\sum_{j=0}^\infty\sum_{k\neq j}\left(3\psi_j^2\psi_{j+\tau_1}^2\psi_k\psi_{k+\tau_2} \right. \\
&\left. +6\psi_j^2\psi_{j+\tau_1}\psi_{j+\tau_2}\psi_k\psi_{k+\tau_1}+3\psi_j\psi_{j+\tau_1}^2\psi_{j+\tau_2}\psi_k^2\right)\E\left[z_t^4z_{t+j-k}^2\right] \\
&+\sum_{j=0}^\infty\sum_{k\neq j}\sum_{l\neq j+\tau_1,l\neq k+\tau_1} \left(6\psi_j^2\psi_k\psi_{k+\tau_1}\psi_l\psi_{l+\tau_2-\tau_1}+3\psi_j^2\psi_k\psi_{k+\tau_2}\psi_l^2\right) \\
&\times\E\left[z_t^2z_{t+j-k}^2z_{t+j-l+\tau_1}^2\right] +\sum_{j=0}^\infty\sum_{k\neq j}\sum_{l\neq j,l\neq k} 6\psi_j\psi_{j+\tau_1}\psi_k\psi_{k+\tau_1}\psi_l\psi_{l+\tau_2} \\
&\times\E\left[z_t^2z_{t+j-k}^2z_{t+j-l}^2\right] \\
&\E\left[x_t^3x_{t+\tau}\right]=\sum_{j=0}^\infty \psi_j^3\psi_{j+\tau}\E\left[z_t^4\right]+\sum_{j=0}^\infty \sum_{k\neq j} 3\psi_j^2\psi_k\psi_{k+\tau}\E\left[z_t^2z_{t+j-k}^2\right]
\end{align*}

\begin{align*}
&\E\left[x_t^3x_{t+\tau_1}x_{t+\tau_2}^2\right]=\sum_{j=0}^\infty \psi_j^3\psi_{j+\tau_1}\psi_{j+\tau_2}^2\E\left[z_t^6\right]+\sum_{j=0}^\infty \sum_{k\neq j+\tau_2} \psi_j^3\psi_{j+\tau_1}\psi_k^2 \\
&\times \E\left[z_t^4z_{t+j-k+\tau_2}^2\right]+\sum_{j=0}^\infty \sum_{k\neq j+\tau_1} 2\psi_j^3\psi_{j+\tau_2}\psi_k\psi_{k+\tau_2-\tau_1}\E\left[z_t^4z_{t+j-k+\tau_1}^2\right]\\
&+\sum_{j=0}^\infty\sum_{k\neq j}\left(3\psi_j^2\psi_{j+\tau_2}^2\psi_k\psi_{k+\tau_1}+6\psi_j^2\psi_{j+\tau_1}\psi_{j+\tau_2}\psi_k\psi_{k+\tau_2} \right. \\
&\left. +3\psi_j\psi_{j+\tau_1}\psi_{j+\tau_2}^2\psi_k^2\right)\E\left[z_t^4z_{t+j-k}^2\right] \\
&+\sum_{j=0}^\infty\sum_{k\neq j}\sum_{l\neq j+\tau_2,l\neq k+\tau_2} 3\psi_j^2\psi_k\psi_{k+\tau_1}\psi_l^2\E\left[z_t^2z_{t+j-k}^2z_{t+j-l+\tau_2}^2\right] \\
&+\sum_{j=0}^\infty\sum_{k\neq j}\sum_{l\neq j+\tau_1,l\neq k+\tau_1} 6\psi_j^2\psi_k\psi_{k+\tau_2}\psi_l\psi_{l+\tau_2-\tau_1}\E\left[z_t^2z_{t+j-k}^2z_{t+j-l+\tau_1}^2\right] \\
&+\sum_{j=0}^\infty\sum_{k\neq j}\sum_{l\neq j,l\neq k} 6\psi_j\psi_{j+\tau_1}\psi_k\psi_{k+\tau_2}\psi_l\psi_{l+\tau_2} \E\left[z_t^2z_{t+j-k}^2z_{t+j-l}^2\right] \\
&\E\left[x_t^2x_{t+\tau}^4\right]=\sum_{j=0}^\infty \psi_j^2\psi_{j+\tau}^4\E\left[z_t^6\right]+\sum_{j=0}^\infty \sum_{k\neq j+\tau} \left(6\psi_j^2\psi_{j+\tau}^2\psi_k^2\E\left[z_t^4z_{t+j-k+\tau}^2\right] \right. \\
&\left. +\psi_j^2\psi_k^4\E\left[z_t^2z_{t+j-k+\tau}^4\right]\right)+\sum_{j=0}^\infty \sum_{k\neq j} 8\psi_j\psi_{j+\tau}^3\psi_k\psi_{k+\tau}\E\left[z_t^4z_{t+j-k}^2\right] \\
&+\sum_{j=0}^\infty\sum_{k\neq j+\tau}\sum_{l\neq j+\tau,l\neq k} 3\psi_j^2\psi_k^2\psi_l^2\E\left[z_t^2z_{t+j-k+\tau}^2z_{t+j-l+\tau}^2\right] \\
&+\sum_{j=0}^\infty\sum_{k\neq j}\sum_{l\neq j,l\neq k} 12\psi_j\psi_{j+\tau}\psi_k\psi_{k+\tau}\psi_l^2\E\left[z_t^2z_{t+j-k}^2z_{t+j-l}^2\right]
\end{align*}

\begin{align*}
&\E\left[x_t^2x_{t+\tau_1}^3x_{t+\tau_2}\right]=\sum_{j=0}^\infty \psi_j^2\psi_{j+\tau_1}^3\psi_{j+\tau_2}\E\left[z_t^6\right]+\sum_{j=0}^\infty \sum_{k\neq j+\tau_1} \left(\left(3\psi_j^2\psi_{j+\tau_1}^2\psi_k\psi_{k+\tau_2-\tau_1} \right.\right. \\
&\left. \left.+3\psi_j^2\psi_{j+\tau_1}\psi_{j+\tau_2}\psi_k^2\right)\E\left[z_t^4z_{t+j-k+\tau_1}^2\right]+\psi_j^2\psi_k^3\psi_{k+\tau_2-\tau_1}\E\left[z_t^2z_{t+j-k+\tau_1}^4\right]\right) \\
&+\sum_{j=0}^\infty\sum_{k\neq j}\left(2\psi_j\psi_{j+\tau_1}^3\psi_k\psi_{k+\tau_2}+6\psi_j\psi_{j+\tau_1}^2\psi_{j+\tau_2}\psi_k\psi_{k+\tau_1}\right)\E\left[z_t^4z_{t+j-k}^2\right] \\
&+\sum_{j=0}^\infty\sum_{k\neq j+\tau_1}\sum_{l\neq j+\tau_1,l\neq k} 3\psi_j^2\psi_k^2\psi_l\psi_{l+\tau_2-\tau_1}\E\left[z_t^2z_{t+j-k+\tau_1}^2z_{t+j-l+\tau_1}^2\right] \\
&+\sum_{j=0}^\infty\sum_{k\neq j}\sum_{l\neq j+\tau_1,l\neq k+\tau_1} \left(6\psi_j\psi_{j+\tau_1}\psi_k\psi_{k+\tau_1}\psi_l\psi_{l+\tau_2-\tau_1}+6\psi_j\psi_{j+\tau_1}\psi_k\psi_{k+\tau_2}\psi_l^2\right) \\
&\times \E\left[z_t^2z_{t+j-k}^2z_{t+j-l+\tau_1}^2\right] \\
&\E\left[x_t^2x_{t+\tau}^2\right]=  \sum_{j=0}^\infty \psi_j^2\psi_{j+\tau}^2\E\left[z_t^4\right]+\sum_{j=0}^\infty \sum_{k\neq j} 2\psi_j\psi_{j+\tau}\psi_k\psi_{k+\tau}\E\left[z_t^2z_{t+j-k}^2\right] \\
&+\sum_{j=0}^\infty \sum_{k\neq j+\tau} \psi_j^2\psi_k^2\E\left[z_t^2z_{t+j-k+\tau}^2\right]  \\
&\E\left[x_t^2x_{t+\tau_1}^2x_{t+\tau_2}^2\right]=\sum_{j=0}^\infty \psi_j^2\psi_{j+\tau_1}^2\psi_{j+\tau_2}^2\E\left[z_t^6\right]+\sum_{j=0}^\infty \sum_{k\neq j+\tau_2} \psi_j^2\psi_{j+\tau_1}^2\psi_k^2 \\
&\times\E\left[z_t^4z_{t+j-k+\tau_2}^2\right] +\sum_{j=0}^\infty \sum_{k\neq j+\tau_1} \left(\left(\psi_j^2\psi_{j+\tau_2}^2\psi_k^2+4\psi_j^2\psi_{j+\tau_1}\psi_{j+\tau_2}\psi_k\psi_{k+\tau_2-\tau_1}\right) \right.\\
&\left.\times\E\left[z_t^4z_{t+j-k+\tau_1}^2\right]+\psi_j^2\psi_k^2\psi_{k+\tau_2-\tau_1}^2\E\left[z_t^2z_{t+j-k+\tau_1}^4\right]\right) \\
&+\sum_{j=0}^\infty \sum_{k\neq j} \left(4\psi_j\psi_{j+\tau_1}^2\psi_{j+\tau_2}\psi_k\psi_{k+\tau_2}+4\psi_j\psi_{j+\tau_1}\psi_{j+\tau_2}^2\psi_k\psi_{k+\tau_1}\right) 
\times\E\left[z_t^4z_{t+j-k}^2\right] \\
&+\sum_{j=0}^\infty\sum_{k\neq j+\tau_1}\sum_{l\neq j+\tau_2,l\neq k+\tau_2-\tau_1} \psi_j^2\psi_k^2\psi_l^2\E\left[z_t^2z_{t+j-k+\tau_1}^2z_{t+j-l+\tau_2}^2\right] \\
&+\sum_{j=0}^\infty\sum_{k\neq j}\sum_{l\neq j+\tau_2,l\neq k+\tau_2} 2\psi_j\psi_{j+\tau_1}\psi_k\psi_{k+\tau_1}\psi_l^2\E\left[z_t^2z_{t+j-k}^2z_{t+j-l+\tau_2}^2\right] \\
&+\sum_{j=0}^\infty\sum_{k\neq j+\tau_1}\sum_{l\neq j+\tau_1,l\neq k} 2\psi_j^2\psi_k\psi_{k+\tau_2-\tau_1}\psi_l\psi_{l+\tau_2-\tau_1} \E\left[z_t^2z_{t+j-k+\tau_1}^2z_{t+j-l+\tau_1}^2\right] \\
&+\sum_{j=0}^\infty\sum_{k\neq j}\sum_{l\neq j+\tau_1,l\neq k+\tau_1}\left( 2\psi_j\psi_{j+\tau_2}\psi_k\psi_{k+\tau_2}\psi_l^2+8\psi_j\psi_{j+\tau_1}\psi_k\psi_{k+\tau_2}\psi_l\psi_{l+\tau_2-\tau_1}\right) \\
&\times\E\left[z_t^2z_{t+j-k+\tau_1}^2z_{t+j-l+\tau_1}^2\right] \\
\end{align*}

\begin{align*}
&\E\left[x_t^2x_{t+\tau_1}^2x_{t+\tau_2}x_{t+\tau_3}\right]=\sum_{j=0}^\infty \psi_j^2\psi_{j+\tau_1}^2\psi_{j+\tau_2}\psi_{j+\tau_3}\E\left[z_t^6\right]+\sum_{j=0}^\infty \sum_{k\neq j+\tau_2} \psi_j^2\psi_{j+\tau_1}^2 \\
&\times \psi_k\psi_{k+\tau_3-\tau_2}\E\left[z_t^4z_{t+j-k+\tau_2}^2\right]+\sum_{j=0}^\infty \sum_{k\neq j+\tau_1}\left(\left(2\psi_j^2\psi_{j+\tau_1}\psi_{j+\tau_2}\psi_k\psi_{k+\tau_3-\tau_1} \right.\right. \\
&\left. +2\psi_j^2\psi_{j+\tau_1} \psi_{j+\tau_3}\psi_k\psi_{k+\tau_2-\tau_1}+\psi_j^2\psi_{j+\tau_2}\psi_{j+\tau_3}\psi_k^2\right)\E\left[z_t^4z_{t+j-k+\tau_1}^2\right]+\psi_j^2\psi_k^2  \\
&\left. \times\psi_{k+\tau_2-\tau_1}\psi_{k+\tau_3-\tau_1}\E\left[z_t^2z_{t+j-k+\tau_1}^4\right]\right)+\sum_{j=0}^\infty \sum_{k\neq j}\left(2\psi_j\psi_{j+\tau_1}^2\psi_{j+\tau_2}\psi_k\psi_{k+\tau_3} \right. \\
&\left. +2\psi_j\psi_{j+\tau_1}^2\psi_{j+\tau_3}\psi_k\psi_{k+\tau_2}+4\psi_j\psi_{j+\tau_1}\psi_{j+\tau_2}\psi_{j+\tau_3}\psi_k\psi_{k+\tau_1}\right)\E\left[z_t^4z_{t+j-k}^2\right] \\
&+\sum_{j=0}^\infty\sum_{k\neq j+\tau_1}\sum_{l\neq j+\tau_2,l\neq k+\tau_2-\tau_1}\psi_j^2\psi_k^2\psi_l\psi_{l+\tau_3-\tau_ 2}\E\left[z_t^2z_{t+j-k+\tau_1}^2z_{t+j-l+\tau_2}^2\right] \\
&+\sum_{j=0}^\infty\sum_{k\neq j}\sum_{l\neq j+\tau_2,l\neq k+\tau_2} 2\psi_j\psi_{j+\tau_1}\psi_k\psi_{k+\tau_1}\psi_l\psi_{l+\tau_3-\tau_2}\E\left[z_t^2z_{t+j-k}^2z_{t+j-l+\tau_2}^2\right] \\
&+\sum_{j=0}^\infty\sum_{k\neq j+\tau_1}\sum_{l\neq j+\tau_1,l\neq k} 2\psi_j^2\psi_k\psi_{k+\tau_2-\tau_1}\psi_l\psi_{l+\tau_3-\tau_1}\E\left[z_t^2z_{t+j-k+\tau_1}^2z_{t+j-l+\tau_1}^2\right] \\
&+\sum_{j=0}^\infty\sum_{k\neq j}\sum_{l\neq j+\tau_1,l\neq k+\tau_1} \left(4\psi_j\psi_{j+\tau_1}\psi_k\psi_{k+\tau_2}\psi_l\psi_{l+\tau_3-\tau_1}+4\psi_j\psi_{j+\tau_1}\psi_k\psi_{k+\tau_3}\psi_l\right. \\
&\left.\times\psi_{l+\tau_2-\tau_1}+2\psi_j\psi_{j+\tau_2}\psi_k\psi_{k+\tau_3}\psi_l^2\right) \E\left[z_t^2z_{t+j-k}^2z_{t+j-l+\tau_1}^2\right] \\
\end{align*}

\begin{align*}
&\E\left[x_t^2x_{t+\tau_1}x_{t+\tau_2}^3\right]=\sum_{j=0}^\infty \psi_j^2\psi_{j+\tau_1}\psi_{j+\tau_2}^3\E\left[z_t^6\right]+\sum_{j=0}^\infty \sum_{k\neq j+\tau_2} 3\psi_j^2\psi_{j+\tau_1}\psi_{j+\tau_2}\psi_k^2 \\
&\times \E\left[z_t^4z_{t+j-k+\tau_2}^2\right]+\sum_{j=0}^\infty \sum_{k\neq j+\tau_1}\left(3\psi_j^2\psi_{j+\tau_2}^2\psi_k\psi_{k+\tau_2-\tau_1}\E\left[z_t^4z_{t+j-k+\tau_1}^2\right]\right. \\
&\left. +\psi_j^2\psi_k\psi_{k+\tau_2-\tau_1}^3\E\left[z_t^2z_{t+j-k+\tau_1}^4\right]\right)+\sum_{j=0}^\infty\sum_{k\neq j}\left(6\psi_j\psi_{j+\tau_1}\psi_{j+\tau_2}^2\psi_k\psi_{k+\tau_2} \right. \\
&\left. +2\psi_j\psi_{j+\tau_2}^3\psi_k\psi_{k+\tau_1}\right)\E\left[z_t^4z_{t+j-k}^2\right] \\
&+\sum_{j=0}^\infty\sum_{k\neq j+\tau_1}\sum_{l\neq j+\tau_2,l\neq k+\tau_2-\tau_1}3\psi_j^2\psi_k\psi_{k+\tau_2-\tau_1}\psi_l^2\E\left[z_t^2z_{t+j-k+\tau_1}^2z_{t+j-l+\tau_2}^2\right] \\
&+\sum_{j=0}^\infty\sum_{k\neq j}\sum_{l\neq j+\tau_2,l\neq k+\tau_2} 6\psi_j\psi_{j+\tau_1}\psi_k\psi_{k+\tau_2}\psi_l^2\E\left[z_t^2z_{t+j-k}^2z_{t+j-l+\tau_2}^2\right] \\
&+\sum_{j=0}^\infty\sum_{k\neq j}\sum_{l\neq j+\tau_1,l\neq k+\tau_1} 6\psi_j\psi_{j+\tau_2}\psi_k\psi_{k+\tau_2}\psi_l\psi_{l+\tau_2-\tau_1} \E\left[z_t^2z_{t+j-k}^2z_{t+j-l+\tau_1}^2\right]
\end{align*}

\begin{align*}
&\E\left[x_t^2x_{t+\tau_1}x_{t+\tau_2}^2x_{t+\tau_3}\right]=\sum_{j=0}^\infty \psi_j^2\psi_{j+\tau_1}\psi_{j+\tau_2}^2\psi_{j+\tau_3}\E\left[z_t^6\right]+\sum_{j=0}^\infty \sum_{k\neq j+\tau_2} \left(2\psi_j^2\psi_{j+\tau_1}\right. \\
&\left.\times \psi_{j+\tau_2} \psi_k\psi_{k+\tau_3-\tau_2}+\psi_j^2\psi_{j+\tau_1}\psi_{j+\tau_3}\psi_k^2\right)\E\left[z_t^4z_{t+j-k+\tau_2}^2\right]+\sum_{j=0}^\infty \sum_{k\neq j+\tau_1}\left(\left(\psi_j^2\psi_{j+\tau_2}^2\right. \right. \\
&\left.  \times \psi_k\psi_{k+\tau_3-\tau_1} +2\psi_j^2\psi_{j+\tau_2}\psi_{j+\tau_3}\psi_k\psi_{k+\tau_2-\tau_1}\right)\E\left[z_t^4z_{t+j-k+\tau_1}^2\right] +\psi_j^2\psi_k\psi_{k+\tau_2-\tau_1}^2 \\
&\left. \times \psi_{k+\tau_3-\tau_1}\E\left[z_t^2z_{t+j-k+\tau_1}^4\right]\right) +\sum_{j=0}^\infty \sum_{k\neq j}\left(2\psi_j\psi_{j+\tau_1}\psi_{j+\tau_2}^2\psi_k\psi_{k+\tau_3}+4\psi_j\psi_{j+\tau_1} \right. \\
&\left. \times \psi_{j+\tau_2}\psi_{j+\tau_3}\psi_k\psi_{k+\tau_2}+2\psi_j\psi_{j+\tau_2}^2\psi_{j+\tau_3}\psi_k\psi_{k+\tau_1}\right)\E\left[z_t^4z_{t+j-k}^2\right] \\
&+\sum_{j=0}^\infty\sum_{k\neq j+\tau_1}\sum_{l\neq j+\tau_2,l\neq k+\tau_2-\tau_1}\left(2\psi_j^2\psi_k\psi_{k+\tau_2-\tau_1}\psi_l\psi_{l+\tau_3-\tau_ 2} +\psi_j^2\psi_k\psi_{k+\tau_3-\tau_1}\psi_l^2\right) \\
&\times \E\left[z_t^2z_{t+j-k+\tau_1}^2z_{t+j-l+\tau_2}^2\right] +\sum_{j=0}^\infty\sum_{k\neq j}\sum_{l\neq j+\tau_2,l\neq k+\tau_2} \left(4\psi_j\psi_{j+\tau_1}\psi_k\psi_{k+\tau_2}\psi_l\psi_{l+\tau_3-\tau_2} \right. \\
&\left. +2\psi_j\psi_{j+\tau_1}\psi_k\psi_{k+\tau_3}\psi_l^2\right)\E\left[z_t^2z_{t+j-k}^2z_{t+j-l+\tau_2}^2\right] +\sum_{j=0}^\infty\sum_{k\neq j}\sum_{l\neq j+\tau_1,l\neq k+\tau_1} \left(2\psi_j\psi_{j+\tau_2} \right. \\
&\left. \times \psi_k\psi_{k+\tau_2}\psi_l\psi_{l+\tau_3-\tau_1}+4\psi_j\psi_{j+\tau_2}\psi_k\psi_{k+\tau_3}\psi_l\psi_{l+\tau_2-\tau_1}\right)\E\left[z_t^2z_{t+j-k}^2z_{t+j-l+\tau_1}^2\right] \\
&\E\left[x_t^2x_{t+\tau_1}x_{t+\tau_2}\right]=\sum_{j=0}^\infty \psi_j^2\psi_{j+\tau_1}\psi_{j+\tau_2}\E\left[z_t^4\right]+\sum_{j=0}^\infty \sum_{k\neq j+\tau_1} \psi_j^2\psi_k\psi_{k+\tau_2-\tau_1} \\
&\times\E\left[z_t^2z_{t+j-k+\tau_1}^2\right]+\sum_{j=0}^\infty \sum_{k\neq j} 2\psi_j\psi_{j+\tau_1}\psi_k\psi_{k+\tau_2}\E\left[z_t^2z_{t+j-k+\tau}^2\right]
\end{align*}

\begin{align*}
&\E\left[x_t^2x_{t+\tau_1}x_{t+\tau_2}x_{t+\tau_3}^2\right]=\sum_{j=0}^\infty \psi_j^2\psi_{j+\tau_1}\psi_{j+\tau_2}\psi_{j+\tau_3}^2\E\left[z_t^6\right]+\sum_{j=0}^\infty \sum_{k\neq j+\tau_3}\psi_j^2\psi_{j+\tau_1}\psi_{j+\tau_2} \\
&\times \psi_k^2\E\left[z_t^4z_{t+j-k+\tau_3}^2\right]+\sum_{j=0}^\infty \sum_{k\neq j+\tau_2} 2\psi_j^2\psi_{j+\tau_1}\psi_{j+\tau_3}\psi_k\psi_{k+\tau_3-\tau_2}\E\left[z_t^4z_{t+j-k+\tau_2}^2\right] \\
&+\sum_{j=0}^\infty \sum_{k\neq j+\tau_1}\left(\left(2\psi_j^2\psi_{j+\tau_2}\psi_{j+\tau_3}\psi_k\psi_{k+\tau_3-\tau_1}+\psi_j^2\psi_{j+\tau_3}^2\psi_k\psi_{k+\tau_2-\tau_1}\right)\E\left[z_t^4z_{t+j-k+\tau_1}^2\right] \right. \\
&\left.+\psi_j^2\psi_k\psi_{k+\tau_2-\tau_1}\psi_{k+\tau_3-\tau_1}^2\E\left[z_t^2z_{t+j-k+\tau_1}^4\right]\right) +\sum_{j=0}^\infty \sum_{k\neq j}\left(4\psi_j\psi_{j+\tau_1}\psi_{j+\tau_2}\psi_{j+\tau_3}\psi_k\psi_{k+\tau_3} \right. \\
&\left.+2\psi_j\psi_{j+\tau_1}\psi_{j+\tau_3}^2\psi_k\psi_{k+\tau_2}+2\psi_j\psi_{j+\tau_2}\psi_{j+\tau_3}^2\psi_k\psi_{k+\tau_1}\right)\E\left[z_t^4z_{t+j-k}^2\right] \\
&+\sum_{j=0}^\infty\sum_{k\neq j+\tau_1}\sum_{l\neq j+\tau_3,l\neq k+\tau_3-\tau_1}\psi_j^2\psi_k\psi_{k+\tau_2-\tau_1}\psi_l^2\E\left[z_t^2z_{t+j-k+\tau_1}^2z_{t+j-l+\tau_3}^2\right] \\
&+\sum_{j=0}^\infty\sum_{k\neq j}\sum_{l\neq j+\tau_3,l\neq k+\tau_3}2\psi_j\psi_{j+\tau_1}\psi_k\psi_{k+\tau_2}\psi_l^2\E\left[z_t^2z_{t+j-k}^2z_{t+j-l+\tau_3}^2\right] \\
&+\sum_{j=0}^\infty\sum_{k\neq j+\tau_1}\sum_{l\neq j+\tau_2,l\neq k+\tau_2-\tau_1} 2\psi_j^2\psi_k\psi_{k+\tau_3-\tau_1}\psi_l\psi_{l+\tau_3-\tau_2}\E\left[z_t^2z_{t+j-k+\tau_1}^2z_{t+j-l+\tau_2}^2\right] \\
&+\sum_{j=0}^\infty\sum_{k\neq j}\sum_{l\neq j+\tau_2,l\neq k+\tau_2} 4\psi_j\psi_{j+\tau_1}\psi_k\psi_{k+\tau_3}\psi_l\psi_{l+\tau_3-\tau_2}\E\left[z_t^2z_{t+j-k}^2z_{t+j-l+\tau_2}^2\right] \\
&+\sum_{j=0}^\infty\sum_{k\neq j}\sum_{l\neq j+\tau_1,l\neq k+\tau_1} \left(4\psi_j\psi_{j+\tau_2}\psi_k\psi_{k+\tau_3}\psi_l\psi_{l+\tau_3-\tau_1}+2\psi_j\psi_{j+\tau_3}\psi_k\psi_{k+\tau_3}\psi_l\psi_{l+\tau_2-\tau_1}\right) \\
&\times \E\left[z_t^2z_{t+j-k}^2z_{t+j-l+\tau_1}^2\right]
\end{align*}

\begin{align*}
&\E\left[x_tx_{t+\tau_1}^4x_{t+\tau_2}\right]=\sum_{j=0}^\infty \psi_j\psi_{j+\tau_1}^4\psi_{j+\tau_2}\E\left[z_t^6\right]+\sum_{j=0}^\infty \sum_{k\neq j+\tau_1} \left(\left( 4\psi_j\psi_{j+\tau_1}^3\psi_k\psi_{k+\tau_2-\tau_1} \right. \right. \\
&\left. \left.+6\psi_j\psi_{j+\tau_1}^2\psi_{j+\tau_2}\psi_k^2\right)\E\left[z_t^4z_{t+j-k+\tau_1}^2\right]+\psi_j\psi_{j+\tau_2}\psi_k +4\psi_j\psi_{j+\tau_2}\psi_k^3\psi_{k+\tau_2-\tau_1}\E\left[z_t^2z_{t+j-k+\tau_1}^4\right]\right) \\
&+\sum_{j=0}^\infty\sum_{k\neq j+\tau_1}\sum_{l\neq j+\tau_1,l\neq k} \left(12\psi_j\psi_{j+\tau_1}\psi_k\psi_{k+\tau_2-\tau_1}\psi_l^2+3\psi_j\psi_{j+\tau_2}\psi_k^2\psi_l^2\right) \times \\
&\ \ \ \E\left[z_t^2z_{t+j-k+\tau_1}^2z_{t+j-l+\tau_1}^2\right] \\
&\E\left[x_tx_{t+\tau}^3\right]=\sum_{j=0}^\infty \psi_j\psi_{j+\tau}^3\E\left[z_t^4\right]+\sum_{j=0}^\infty \sum_{k\neq j+\tau} 3\psi_j\psi_{j+\tau}\psi_k^2\E\left[z_t^2z_{t+j-k+\tau}^2\right] \\
&\E\left[x_tx_{t+\tau_1}^3x_{t+\tau_2}^2\right]=\sum_{j=0}^\infty \psi_j\psi_{j+\tau_1}^3\psi_{j+\tau_2}^2\E\left[z_t^6\right]+\sum_{j=0}^\infty \sum_{k\neq j+\tau_2} \psi_j\psi_{j+\tau_1}^3\psi_k^2\E\left[z_t^4z_{t+j-k+\tau_2}^2\right] \\
&+\sum_{j=0}^\infty \sum_{k\neq j+\tau_1} \left(\left(6\psi_j\psi_{j+\tau_1}^2\psi_{j+\tau_2}\psi_k\psi_{k+\tau_2-\tau_1} +3\psi_j\psi_{j+\tau_1}\psi_{j+\tau_2}^2\psi_k^2\right)\E\left[z_t^4z_{t+j-k+\tau_1}^2\right] \right. \\
&\left. +\left(3\psi_j\psi_{j+\tau_1}\psi_k^2\psi_{k+\tau_2-\tau_1}^2+2\psi_j\psi_{j+\tau_2}\psi_k^3\psi_{k+\tau_2-\tau_1}\right)\E\left[z_t^2z_{t+j-k+\tau_1}^4\right]\right) \\
&+\sum_{j=0}^\infty\sum_{k\neq j+\tau_1}\sum_{l\neq j+\tau_2,l\neq k+\tau_2-\tau_1} 3\psi_j\psi_{j+\tau_1}\psi_k^2\psi_l^2\E\left[z_t^2z_{t+j-k+\tau_1}^2z_{t+j-l+\tau_2}^2\right] \\
&+\sum_{j=0}^\infty\sum_{k\neq j+\tau_1}\sum_{l\neq j+\tau_1,l\neq k} \left(6\psi_j\psi_{j+\tau_1}\psi_k\psi_{k+\tau_2-\tau_1}\psi_l\psi_{l+\tau_2-\tau_1}+6\psi_j\psi_{j+\tau_2}\psi_k^2\psi_l\psi_{l+\tau_2-\tau_1}\right) \\
&\times \E\left[z_t^2z_{t+j-k}^2z_{t+j-l+\tau_1}^2\right]
\end{align*}

\begin{align*}
&\E\left[x_tx_{t+\tau_1}^2x_{t+\tau_2}^3\right]=\sum_{j=0}^\infty \psi_j\psi_{j+\tau_1}^2\psi_{j+\tau_2}^3\E\left[z_t^6\right]+\sum_{j=0}^\infty \sum_{k\neq j+\tau_2} 3\psi_j\psi_{j+\tau_1}^2\psi_{j+\tau_2}\psi_k^2 \\
&\times \E\left[z_t^4z_{t+j-k+\tau_2}^2\right]+\sum_{j=0}^\infty \sum_{k\neq j+\tau_1}\left(\left(6\psi_j\psi_{j+\tau_1}\psi_{j+\tau_2}^2\psi_k\psi_{k+\tau_2-\tau_1}+\psi_j\psi_{j+\tau_2}^3\psi_k^2\right) \right. \\
& \times \E\left[z_t^4z_{t+j-k+\tau_1}^2\right]+\left(2\psi_j\psi_{j+\tau_1}\psi_k\psi_{k+\tau_2-\tau_1}^3+3\psi_j\psi_{j+\tau_2}\psi_k^2\psi_{k+\tau_2-\tau_1}^2\right) \\
&\left.\times \E\left[z_t^2z_{t+j-k+\tau_1}^4\right]\right)+\sum_{j=0}^\infty\sum_{k\neq j+\tau_1}\sum_{l\neq j+\tau_2,l\neq k+\tau_2-\tau_1} \left(6\psi_j\psi_{j+\tau_1}\psi_k\psi_{k+\tau_2-\tau_1}\psi_l^2 \right. \\
&\left.+3\psi_j\psi_{j+\tau_2}\psi_k^2\psi_l^2\right)\E\left[z_t^2z_{t+j-k+\tau_1}^2z_{t+j-l+\tau_2}^2\right] \\
&+\sum_{j=0}^\infty\sum_{k\neq j+\tau_1}\sum_{l\neq j+\tau_1,l\neq k} 6\psi_j\psi_{j+\tau_2}\psi_k\psi_{k+\tau_2-\tau_1}\psi_l\psi_{l+\tau_2-\tau_1}\E\left[z_t^2z_{t+j-k+\tau_1}^2z_{t+j-l+\tau_1}^2\right] \\
&\E\left[x_tx_{t+\tau_1}^2x_{t+\tau_2}^2x_{t+\tau_3}\right]=\sum_{j=0}^\infty \psi_j\psi_{j+\tau_1}^2\psi_{j+\tau_2}^2\psi_{j+\tau_3}\E\left[z_t^6\right]+\sum_{j=0}^\infty \sum_{k\neq j+\tau_2}\left(2\psi_j\psi_{j+\tau_1}^2 \right. \\
&\left. \ \ \ \psi_{j+\tau_2}\psi_k\psi_{k+\tau_3-\tau_2}+\psi_j\psi_{j+\tau_1}^2\psi_{j+\tau_3}\psi_k^2\right)\E\left[z_t^4z_{t+j-k+\tau_2}^2\right] \\
&+\sum_{j=0}^\infty \sum_{k\neq j+\tau_1}\left(\left(2\psi_j\psi_{j+\tau_1}\psi_{j+\tau_2}^2\psi_k\psi_{k+\tau_3-\tau_1}+4\psi_j\psi_{j+\tau_1}\psi_{j+\tau_2}\psi_{j+\tau_3}\psi_k\psi_{k+\tau_2-\tau_1} \right.\right. \\
&\left. +\psi_j\psi_{j+\tau_2}^2\psi_{j+\tau_3}\psi_k^2\right)\E\left[z_t^4z_{t+j-k+\tau_1}^2\right]+\left(2\psi_j\psi_{j+\tau_1}\psi_k\psi_{k+\tau_2-\tau_1}^2\psi_{k+\tau_3-\tau_1} \right. \\
&\left. \left. +2\psi_j\psi_{j+\tau_2}\psi_k^2\psi_{k+\tau_2-\tau_1}\psi_{k+\tau_3-\tau_1}+\psi_j\psi_{j+\tau_3}\psi_k^2\psi_{k+\tau_2-\tau_1}^2\right)\E\left[z_t^2z_{t+j-k+\tau_1}^4\right]\right) \\
&+\sum_{j=0}^\infty\sum_{k\neq j+\tau_1}\sum_{l\neq j+\tau_2,l\neq k+\tau_2-\tau_1}\left(4\psi_j\psi_{j+\tau_1}\psi_k\psi_{k+\tau_2-\tau_1}\psi_l\psi_{l+\tau_3-\tau_2}+2\psi_j\psi_{j+\tau_1} \right. \\
& \left. \times \psi_k\psi_{k+\tau_3-\tau_1}\psi_l^2+2\psi_j\psi_{j+\tau_2}\psi_k^2\psi_l\psi_{l+\tau_3-\tau_2}+\psi_j\psi_{j+\tau_3}\psi_k^2\psi_l^2\right) \\
&\times \E\left[z_t^2z_{t+j-k+\tau_1}^2z_{t+j-l+\tau_2}^2\right]+\sum_{j=0}^\infty\sum_{k\neq j+\tau_1}\sum_{l\neq j+\tau_1,l\neq k} \left(4\psi_j\psi_{j+\tau_2} \psi_k\psi_{k+\tau_2-\tau_1}  \right. \\
&\left. \times \psi_l\psi_{l+\tau_3-\tau_1}+2\psi_j\psi_{j+\tau_3}\psi_k\psi_{k+\tau_2-\tau_1}\psi_l\psi_{l+\tau_2-\tau_1}\right)\E\left[z_t^2z_{t+j-k+\tau_1}^2z_{t+j-l+\tau_1}^2\right]
\end{align*}

\begin{align*}
&\E\left[x_tx_{t+\tau_1}^2x_{t+\tau_2}\right]=\sum_{j=0}^\infty \psi_j\psi_{j+\tau_1}^2\psi_{j+\tau_2}\E\left[z_t^4\right]+\sum_{j=0}^\infty \sum_{k\neq j+\tau_1}\left(2\psi_j\psi_{j+\tau_1}\times \right. \\
&\left. \ \ \ \psi_k\psi_{k+\tau_2-\tau_1}+2\psi_j\psi_{j+\tau_2}\psi_k^2\right)\E\left[z_t^2z_{t+j-k+\tau_1}^2\right] \\
&\E\left[x_tx_{t+\tau_1}^2x_{t+\tau_2}x_{t+\tau_3}^2\right]=\sum_{j=0}^\infty \psi_j\psi_{j+\tau_1}^2\psi_{j+\tau_2}\psi_{j+\tau_3}^2\E\left[z_t^6\right]+\sum_{j=0}^\infty \sum_{k\neq j+\tau_3}\psi_j\psi_{j+\tau_1}^2\psi_{j+\tau_2} \\
&\times \psi_k^2\E\left[z_t^4z_{t+j-k+\tau_3}^2\right]+\sum_{j=0}^\infty \sum_{k\neq j+\tau_2}2\psi_j\psi_{j+\tau_1}^2\psi_{j+\tau_3}\psi_k\psi_{k+\tau_3-\tau_2}\E\left[z_t^4z_{t+j-k+\tau_2}^2\right] \\
&+\sum_{j=0}^\infty \sum_{k\neq j+\tau_1}\left(\left(4\psi_j\psi_{j+\tau_1}\psi_{j+\tau_2}\psi_{j+\tau_3}\psi_k\psi_{k+\tau_3-\tau_1}+2\psi_j\psi_{j+\tau_1}\psi_{j+\tau_3}^2\psi_k\psi_{k+\tau_2-\tau_1} \right.\right. \\
&\left. +\psi_j\psi_{j+\tau_2}\psi_{j+\tau_3}^2\psi_k^2\right)\E\left[z_t^4z_{t+j-k+\tau_1}^2\right]+\left(2\psi_j\psi_{j+\tau_1}\psi_k\psi_{k+\tau_2-\tau_1}\psi_{k+\tau_3-\tau_1}^2 \right. \\
&\left. \left. +\psi_j\psi_{j+\tau_2}\psi_k^2\psi_{k+\tau_3-\tau_1}^2+2\psi_j\psi_{j+\tau_3}\psi_k^2\psi_{k+\tau_2-\tau_1}\psi_{k+\tau_3-\tau_1}\right)\E\left[z_t^2z_{t+j-k+\tau_1}^4\right]\right) \\
&+\sum_{j=0}^\infty\sum_{k\neq j+\tau_1}\sum_{l\neq j+\tau_3,l\neq k+\tau_3-\tau_1}\left(2\psi_j\psi_{j+\tau_1}\psi_k\psi_{k+\tau_2-\tau_1}\psi_l^2+\psi_j\psi_{j+\tau_2}\psi_k^2\psi_l^2\right) \\
&\times \E\left[z_t^2z_{t+j-k+\tau_1}^2z_{t+j-l+\tau_2}^2\right]+\sum_{j=0}^\infty\sum_{k\neq j+\tau_1}\sum_{l\neq j+\tau_2,l\neq k+\tau_2-\tau_1}\left(4\psi_j\psi_{j+\tau_1}\psi_k\psi_{k+\tau_3-\tau_1} \right. \\
&\left. \times \psi_l\psi_{l+\tau_3-\tau_2}+2\psi_j\psi_{j+\tau_3}\psi_k^2\psi_l\psi_{l+\tau_3-\tau_2}\right)\E\left[z_t^2z_{t+j-k+\tau_1}^2z_{t+j-l+\tau_2}^2\right] \\
&+\sum_{j=0}^\infty\sum_{k\neq j+\tau_1}\sum_{l\neq j+\tau_1,l\neq k} \left(2\psi_j\psi_{j+\tau_2}\psi_k\psi_{k+\tau_3-\tau_1}\psi_l\psi_{l+\tau_3-\tau_1}+4\psi_j\psi_{j+\tau_3}\psi_k\psi_{k+\tau_2-\tau_1} \right. \\
&\left. \times \psi_l\psi_{l+\tau_3-\tau_1}\right)\E\left[z_t^2z_{t+j-k+\tau_1}^2z_{t+j-l+\tau_1}^2\right] \\
&\E\left[x_tx_{t+\tau}\right]=\sum_{j=0}^\infty \psi_j\psi_{j+\tau}
\end{align*}

\begin{align*}
&\E\left[x_tx_{t+\tau_1}x_{t+\tau_2}^4\right]=\sum_{j=0}^\infty \psi_j\psi_{j+\tau_1}\psi_{j+\tau_2}^4\E\left[z_t^6\right]+\sum_{j=0}^\infty \sum_{k\neq j+\tau_2} \left(6\psi_j\psi_{j+\tau_1}\psi_{j+\tau_2}^2\psi_k^2 \right. \\
&\left. \times \E\left[z_t^4z_{t+j-k+\tau_2}^2\right]+\psi_j\psi_{j+\tau_1}\psi_k^4\E\left[z_t^2z_{t+j-k+\tau_2}^4\right]\right)+\sum_{j=0}^\infty \sum_{k\neq j+\tau_1}\left(4\psi_j\psi_{j+\tau_2}^3\psi_k \right. \\
&\left. \times \psi_{k+\tau_2-\tau_1}\E\left[z_t^4z_{t+j-k+\tau_1}^2\right]+4\psi_j\psi_{j+\tau_2}\psi_k\psi_{k+\tau_2-\tau_1}^3\E\left[z_t^2z_{t+j-k+\tau_1}^4\right]\right) \\
&+\sum_{j=0}^\infty\sum_{k\neq j+\tau_2}\sum_{l\neq j+\tau_2,l\neq k} 3\psi_j\psi_{j+\tau_1}\psi_k^2\psi_l^2\E\left[z_t^2z_{t+j-k+\tau_2}^2z_{t+j-l+\tau_2}^2\right] \\
&+\sum_{j=0}^\infty\sum_{k\neq j+\tau_1}\sum_{l\neq j+\tau_2,l\neq k+\tau_2-\tau_1} 12\psi_j\psi_{j+\tau_2}\psi_k\psi_{k+\tau_2-\tau_1}\psi_l^2\E\left[z_t^2z_{t+j-k+\tau_1}^2z_{t+j-l+\tau_2}^2\right] \\
&\E\left[x_tx_{t+\tau_1}x_{t+\tau_2}^2\right]=\sum_{j=0}^\infty \psi_j\psi_{j+\tau_1}\psi_{j+\tau_2}^2\E\left[z_t^4\right]+\sum_{j=0}^\infty \sum_{k\neq j+\tau_2}\psi_j\psi_{j+\tau_1}\psi_k^2 \\
&\times \E\left[z_t^2z_{t+j-k+\tau_2}^2\right]+\sum_{j=0}^\infty \sum_{k\neq j+\tau_1}2\psi_j\psi_{j+\tau_2}\psi_k\psi_{k+\tau_2-\tau_1}\E\left[z_t^2z_{t+j-k+\tau_1}^2\right]
\end{align*}

\begin{align*}
&\E\left[x_tx_{t+\tau_1}x_{t+\tau_2}^2x_{t+\tau_3}^2\right]=\sum_{j=0}^\infty \psi_j\psi_{j+\tau_1}\psi_{j+\tau_2}^2\psi_{j+\tau_3}^2\E\left[z_t^6\right]+\sum_{j=0}^\infty \sum_{k\neq j+\tau_3}\psi_j\psi_{j+\tau_1} \\
&\times \psi_{j+\tau_2}^2\psi_k^2\E\left[z_t^4z_{t+j-k+\tau_3}^2\right]+\sum_{j=0}^\infty \sum_{k\neq j+\tau_2}\left(\left(4\psi_j\psi_{j+\tau_1}\psi_{j+\tau_2}\psi_{j+\tau_3}\psi_k\psi_{k+\tau_3-\tau_2} \right.\right. \\
&\left. \left. +\psi_j\psi_{j+\tau_1}\psi_{j+\tau_3}^2\psi_k^2\right)\E\left[z_t^4z_{t+j-k+\tau_2}^2\right]+\psi_j\psi_{j+\tau_1}\psi_k^2\psi_{k+\tau_3-\tau_2}^2\E\left[z_t^2z_{t+j-k+\tau_2}^4\right]\right) \\
&+\sum_{j=0}^\infty \sum_{k\neq j+\tau_1}\left(\left(2\psi_j\psi_{j+\tau_2}^2\psi_{j+\tau_3}\psi_k\psi_{k+\tau_3-\tau_1}+2\psi_j\psi_{j+\tau_2}\psi_{j+\tau_3}^2\psi_k\psi_{k+\tau_2-\tau_1}\right) \right. \\
& \times\E\left[z_t^4z_{t+j-k+\tau_1}^2\right]+\left(2\psi_j\psi_{j+\tau_2}\psi_k\psi_{k+\tau_2-\tau_1}\psi_{k+\tau_3-\tau_1}^2+2\psi_j\psi_{j+\tau_3}\psi_k\psi_{k+\tau_2-\tau_1}^2 \right. \\
&\left. \left.\times \psi_{k+\tau_3-\tau_1}\right)\E\left[z_t^2z_{t+j-k+\tau_1}^4\right]\right)+\sum_{j=0}^\infty\sum_{k\neq j+\tau_2}\sum_{l\neq j+\tau_3,l\neq k+\tau_3-\tau_2}\psi_j\psi_{j+\tau_1}\psi_k^2\psi_l^2 \\
&\times \E\left[z_t^2z_{t+j-k+\tau_2}^2z_{t+j-l+\tau_3}^2\right]+\sum_{j=0}^\infty\sum_{k\neq j+\tau_1}\sum_{l\neq j+\tau_3,l\neq k+\tau_3-\tau_1}2\psi_j\psi_{j+\tau_2}\psi_k\psi_{k+\tau_2-\tau_1} \\
&\times \psi_l^2\E\left[z_t^2z_{t+j-k+\tau_1}^2z_{t+j-l+\tau_3}^2\right]+\sum_{j=0}^\infty\sum_{k\neq j+\tau_2}\sum_{l\neq j+\tau_2,l\neq k}2\psi_j\psi_{j+\tau_1}\psi_k\psi_{k+\tau_3-\tau_2}\psi_l \ \\
&\times \psi_{l+\tau_3-\tau_2}\E\left[z_t^2z_{t+j-k+\tau_2}^2z_{t+j-l+\tau_2}^2\right] +\sum_{j=0}^\infty\sum_{k\neq j+\tau_1}\sum_{l\neq j+\tau_2,l\neq k+\tau_2-\tau_1} \left(4\psi_j\psi_{j+\tau_2} \right. \\
&\times \psi_k\psi_{k+\tau_3-\tau_1}\psi_l\psi_{l+\tau_3-\tau_2}+4\psi_j\psi_{j+\tau_3}\psi_k\psi_{k+\tau_2-\tau_1}\psi_l\psi_{l+\tau_3-\tau_2}+2\psi_j\psi_{j+\tau_3}\psi_k \\
&\left. \times \psi_{k+\tau_3-\tau_1}\psi_l^2\right)\E\left[z_t^2z_{t+j-k+\tau_1}^2z_{t+j-l+\tau_2}^2\right]
\end{align*}

\subsection*{Proof of Theorem 4}

(i): Write $\bo V=(\bo v_1,\dots,\bo v_p)^T$ and $\bo v_i=(v_{i1},\dots,v_{ip})^T$. By using the basic properties of covariance, we have
\begin{align*}
\bo G_m(\bo V\bo x)=&\sum_{i=1}^p\sum_{j=1}^p\cov^2(\bo V\bo x_t\bo x_t^T\bo V^T,(\bo V\bo x)_{t-l,i}(\bo V\bo x)_{t-l,j}) \\
=&\sum_{i=1}^p\sum_{j=1}^p\cov^2(\bo V\bo x_t\bo x_t^T\bo V^T,\bo v_i^T\bo x_{t-l}\bo v_j^T\bo x_{t-l}) \\
=&\sum_{i=1}^p\sum_{j=1}^p\cov(\bo V\bo x_t\bo x_t^T\bo V^T,\sum_{i'=1}^p (v_{ii'}x_{t-l,i'})\sum_{j'=1}^p (v_{jj'}x_{t-l,j'})) \\
&\cov(\bo V\bo x_t\bo x_t^T\bo V^T,\sum_{i''=1}^p (v_{ii''}x_{t-l,i''})\sum_{j''=1}^p (v_{jj''}x_{t-l,j''}))^T \\
=&\sum_{i=1}^p\sum_{j=1}^p\sum_{i'=1}^p\sum_{j'=1}^p\sum_{i''=1}^p\sum_{j''=1}^p\cov(\bo V\bo x_t\bo x_t^T\bo V^T, v_{ii'}v_{jj'}x_{t-l,i'}x_{t-l,j'}) \\
&\cov(\bo V\bo x_t\bo x_t^T\bo V^T, v_{ii''}v_{jj''}x_{t-l,i''}x_{t-l,j''})^T \\
=& \sum_{i=1}^p\sum_{j=1}^p\sum_{i'=1}^p\sum_{j'=1}^p\sum_{i''=1}^p\sum_{j''=1}^p v_{ii'}v_{jj'}v_{ii''}v_{jj''}
\cov(\bo V\bo x_t\bo x_t^T\bo V^T, x_{t-l,i'}x_{t-l,j'}) \\
&\cov(\bo V\bo x_t\bo x_t^T\bo V^T, x_{t-l,i''}x_{t-l,j''})^T \\
=& \sum_{i'=1}^p\sum_{j'=1}^p\sum_{i''=1}^p\sum_{j''=1}^p \cov(\bo V\bo x_t\bo x_t^T\bo V^T, x_{t-l,i'}x_{t-l,j'})
\cov(\bo V\bo x_t\bo x_t^T\bo V^T, x_{t-l,i''}x_{t-l,j''})^T \\
&\sum_{i=1}^p(v_{ii'}v_{ii''})\sum_{j=1}^p(v_{jj'}v_{jj''}).
\end{align*}
Since $\bo V$ is orthogonal, $\sum_{i=1}^p(v_{ii'}v_{ii''})=1$, if $i'=i''$, and 0 otherwise. Hence,
\begin{align*}
\bo G_m(\bo V\bo x)=& \sum_{i'=1}^p\sum_{j'=1}^p \cov(\bo V\bo x_t\bo x_t^T\bo V^T, x_{t-l,i'}x_{t-l,j'}) \cov(\bo V\bo x_t\bo x_t^T\bo V^T, x_{t-l,i'}x_{t-l,j'}) \\
=&\sum_{i'=1}^p\sum_{j'=1}^p \bo V\cov(\bo x_t\bo x_t^T, x_{t-l,i'}x_{t-l,j'})\bo V^T \bo V\cov(\bo x_t\bo x_t^T, x_{t-l,i'}x_{t-l,j'})^T\bo V^T \\
=&\bo V\left(\sum_{i'=1}^p\sum_{j'=1}^p \cov^2(\bo x_t\bo x_t^T, x_{t-l,i'}x_{t-l,j'})\right)\bo V^T \\
=&\bo V\bo G_m(\bo x)\bo V^T.
\end{align*}

(ii): Let us write
\[
\bo G_m(\bo x)=\sum_{l=1}^m \sum_{i=1}^p\sum_{j=1}^p \bo g_l^{ij}(\bo x)\bo g_l^{ij}(\bo x)^T,
\]
where $\bo g_l^{ij}(\bo x)=\cov(\bo x_t\bo x_t^T,x_{t-l,i}x_{t-l,j})$. It is easy to see that
\[
\bo g_l^{ii}(\bo s)_{jk}=\cov(\bo s_{tj}\bo s_{tk}^T,s_{t-l,i}s_{t-l,i})=0,
\]
if $j\neq k$, and thus $\bo g_l^{ii}(\bo s)$ and $\bo g_l^{ii}(\bo s)\bo g_l^{ii}(\bo s)^T$ are diagonal.
For $i\neq r$,
\[
\bo g_l^{ir}(\bo s)_{jk}=\cov(\bo s_{tj}\bo s_{tk}^T,s_{t-l,i}s_{t-l,r})
\]
can be nonzero only if $j=i$ and $k=r$, or $j=r$ and $k=i$. This implies that the only possibly nonzero elements of
$\bo g_l^{ir}(\bo s)\bo g_l^{ir}(\bo s)^T$ are the $i$th and $r$th diagonal elements. Therefore also $\bo G_m(\bo s)$
is diagonal as a sum of diagonal matrices.

\end{document}